\documentclass[review,onefignum,onetabnum,nolineno]{siamart220329}
\nolinenumbers
\usepackage[utf8]{inputenc}
\usepackage{amsfonts}
\usepackage{amsmath, bm}
\usepackage{graphicx}
\usepackage{color}
\usepackage{tikz}
\usepackage{hyperref}
\usepackage{amssymb}
\usepackage{xfrac}
\usepackage{babel}
\usepackage[font=small,labelfont=bf]{caption}

\newcommand{\xx}{\bm{x}}
\newcommand{\uu}{\bm{u}}
\newcommand{\ff}{\bm{f}}

\newcommand{\nn}{\bm{n}}

\newcommand{\ZZ}{{\mathbb Z}}
\newcommand{\NN}{{\mathbb N}}

\newcommand{\corange}[1]{{\color{blue}  #1}}

\newsiamremark{remark}{Remark} 

\title{Dense cell-by-cell systems of PDEs:
approximation, spectral analysis, and preconditioning\thanks{Submitted to the editors
September 2024.}}

\author{Pietro Benedusi\thanks{Euler Institute, Universit\`a della Svizzera Italiana, Lugano, Switzerland
(\email{benedp@usi.ch}).}
\and Paola Ferrari\thanks{Department of Mathematics, University of Genoa, Genoa, Italy\\
School of Mathematics and Natural Sciences, University of Wuppertal, Wuppertal, Germany}
\and Marius Causemann\thanks{Simula Research Laboratory, Oslo, Norway} 
\and Stefano Serra-Capizzano\thanks{University of Insubria, Como, Italy \\
Uppsala University, Uppsala, Sweden} 
}

\begin{document}
\maketitle
	
\begin{abstract}
In the present study, we consider the Extra-Membrane-Intra model (EMI) for the simulation of excitable tissues at the cellular level. We provide the (possibly large) system of partial differential equations (PDEs), equipped with ad hoc boundary conditions, relevant to model portions of excitable tissues, composed of several cells. In particular, we study two geometrical settings: computational cardiology and neuroscience.
The Galerkin approximations to the considered system of PDEs lead to large linear systems of algebraic equations, where the coefficient  matrices depend on the number $N$ of cells and the fineness parameters. We give a structural and spectral analysis of the related matrix-sequences with $N$ fixed and with fineness parameters tending to zero. Based on the theoretical results, we propose preconditioners and specific multilevel solvers.
Numerical experiments are presented and critically discussed,  showing that a monolithic multilevel solver is efficient and robust with respect to all the problem and discretization parameters. In particular, we include numerical results increasing the number of cells $N$, both for idealized geometries (with $N$ exceeding $10^5$) and for realistic, densely populated 3D tissue reconstruction.

\end{abstract}


\section{Introduction}

The cell-by-cell framework, also known as EMI (Extra, Intra, Membrane) model, is a successful tool in computational electrophysiology, employed to simulate excitable tissues at a cellular scale.
Compared to homogenized models, such as the well-known mono- and bidomain equations, its essential novelty consists in resolving cell morphologies explicitly, enabling detailed biological simulations. In contrast, homogenized models assume extra- and intracellular quantities to be defined in the whole spatial domain without a geometrical representation of cell membranes. In this sense, EMI models allow a more detailed modelling of cellular interplay, for example resolving complex cell networks or capturing inhomogeneities along cells membranes. 
The interested reader is referred to \cite{EMI_book}, where an exhaustive description of the EMI model, its derivation, some theoretical results, and several applications are presented and discussed in detail.
This framework is especially well suited to simulate cellular activity in the human brain, where the cell-scale morphology is highly complex \cite{benedusi2024scalable,ellingsrud2024splitting}. In this context, 3D dense reconstructions of neural tissue have recently become available \cite{velicky2023dense,motta2019dense, shapson2024petavoxel,cali20193d,bartol2015computational}. Moreover, in the field of computational cardiology, we see a growing interest in modeling cardiomyocytes, with their structured pattern \cite{rosilhoBEM,huynh2022convergence,jaeger2021efficient,jaeger2021millimeters,tveito2017cell}.
In both applications, it is essential to accurately resolve the propagation of electrical excitation in the tissue under examination \cite{buccino2021improving,ellingsrud2021cell,ellingsrud2020finite,de2023boundary,mori2009numerical}.

Considering the realistic setting of a tissue composed of several cells, with possibly different properties, we provide the associated system of partial differential equations (PDEs) coupled via membrane ordinary differential equations ODEs, and equipped with ad hoc boundary conditions. 
Subsequently, we consider Galerkin-type approxi\-mations of the PDEs system, leading to large linear systems of algebraic equations, where the resulting coefficient matrices depend on the number $N$ of cells and on the spatial discretization. 
The resulting matrices present different levels of structure: at the outermost level, we have a block structure (e.g. a block arrowhead structure in the context of neuroscience applications), while each block has a (block) two-level Toeplitz \cite{tyrt,TilliNota} or (block) two-level Generalized Locally Toeplitz (GLT) nature \cite{GS-II,block-glt-dD}. The number of levels is dictated by the dimensionality of the domain. For example, in the case of three-dimensional domains a (block) three-level Toeplitz or (block) three-level GLT (asymptotical) structure would occur in any single block.

By generalizing the work in \cite{benedusi2023EMI}, concerning the case of $N=1$, our study is developed in the following two directions:
 \begin{itemize}
\item[(i)] first we give a structural and spectral analysis of the EMI matrix-sequences with $N$ fixed and with fineness parameters tending to zero;
\item[(ii)] based on the the theoretical results, we propose preconditioners and specific multilevel solvers, tested in different application settings, both for idealized and realistic geometries, varying the number of cells $N$ and discretization parameters.
\end{itemize}
We emphasize that the GLT analysis is not completely sufficient for the considered setting and hence specific theoretical tools are developed here, allowing us to give a quite complete picture of the spectral features  of the resulting matrices and matrix-sequences (e.g. the eigenvalue distribution). Numerical experiments are presented and critically discussed, showing the efficiency and robustness of a proposed monolithic multilevel strategy. As a benchmark, we highlight \cite{jaeger2021efficient}, where an efficient solution strategy is described for an idealized cardiac geometry, discretized with finite differen\-ces, increasing the number of cells up to $N>10^5$.

The present work is organized as follows. In Section \ref{sec:EMI-pb} the continuous problem is introduced together with possible generalizations.
Section \ref{sec:approx} considers a basic Galerkin strategy for approximating the examined problem. The spectral analysis is given in Section \ref{sec:spectral} in terms of distribution results and degenerating eigenspaces. Specifically, Subsection \ref{ssec:GLTbackground} lays out the foundational theories and concepts necessary to understand the distribution of the entire system (presented in a general form in Subsection \ref{ssec:symbol}) and the specific stiffness matrices and matrix-sequences (discussed in Subsection \ref{ssec:symbol-bis}). These findings represent the starting point for proposing in Section \ref{sec:num} a preconditioned multilevel strategy, where a selection of numerical experiments is discussed, in connection with the spectral analysis and the related algorithmic proposals.
Finally, Section \ref{sec:final} contains conclusions and a list of relevant open problems.

\section{The EMI model}\label{sec:EMI-pb}

We introduce the partial differential equations characteri\-zing the EMI model. Compared to homogenized models, we observe that the novelty of the EMI approach lies in the fact that the cellular membrane $\Gamma$ is explicitly represented, as well as the intra- and extra-cellular subdomains.

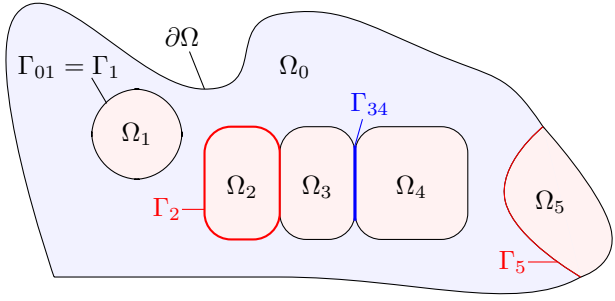
\begin{figure}
    \centering
    \begin{tikzpicture}

      \draw[fill=blue!5]  plot[very thick, smooth, tension=.9] coordinates {(-6,0) (-6.5,3.5) (-4,2.5) (-3,3.5) (-1,3) (0.5,2) (1,0) (-3,0) (-6,0)};

  \draw[rounded corners=18pt, fill=red!5]
  (-5.5,1.3) rectangle ++(1.2,1.2);

 \draw[rounded corners=10pt,fill=red!5]
      (-3,0.5) rectangle ++(1,1.5);

\draw[rounded corners=10pt,fill=red!5]
      (-2,0.5) rectangle ++(1.5,1.5);

\draw[rounded corners=10pt, red, fill=red!5,thick]
  (-4,0.5) rectangle ++(1,1.5);

      \draw[fill=red!5]  plot[very thick, smooth, tension=.9] coordinates { (0.5,2) (0,1) (1,0)};

      \draw[fill=red!5]  plot[very thick, smooth, tension=.9] coordinates { (0.5,2) (1.4,0.65)  (1,0)};

      \draw[color=red]  plot[very thick, smooth, tension=.9] coordinates { (0.5,2) (0,1) (1,0)};

    \node at (-2.8,2.8) {$\Omega_0$};
    \node at (-4.9,1.9) {$\Omega_1$};
    \node at (-3.5,1.2) {$\Omega_2$};
    \node at (-2.5,1.2) {$\Omega_3$};
    \node at (-1.25,1.2) {$\Omega_4$};
    \node at (0.6,1) {$\Omega_5$};
    \node at (-4.3,3.2) {$\partial\Omega$};
    \node[red] at (-4.5,0.9) {$\Gamma_2$};
    \node at (-5.8,2.8) {$\Gamma_{01}=\Gamma_1$};
    \node[red] at (0.1,0.2) {$\Gamma_5$};
    \node[blue] at (-1.8,2.3) {$\Gamma_{34}$};

    \draw[] (-4,2.5) -- (-4.2,3);
    \draw[red] (-4,0.9) -- (-4.3,0.9);
    \draw[] (-5.3,2.31) -- (-5.5,2.65);
    \draw[red] (0.3,0.2) -- (0.69,0.2);
    \draw[blue] (-2,1.75) -- (-1.9,2.1);

    \draw[blue,very thick] (-2,0.75) -- (-2,1.75);
    
    \end{tikzpicture}
    \caption{Example of an EMI geometry with $N=5$ cells with $\Omega=\left(\bigcup_{i=0}^5\Omega_i \right)\cup\left(\bigcup_{i=1}^5\Gamma_i \right)$. Here, boundary conditions are enforced only for $\Omega_0$ and $\Omega_5$ and $\Gamma_{34}$ is an example of a common interface between two cells (e.g. a gap junction).}
    \label{fig:EMI}
\end{figure}

We consider a domain $\Omega\subset\mathbb{R}^d$, with typically $d\in\{2,3\}$, composed of $N$ disjoint domains, or cells, $\Omega_1,\Omega_2,\ldots,\Omega_N$,  an extracellular media $\Omega_0$ surrounding the $N$ cells, so that  $\overline\Omega=\bigcup_{i=0}^{N}\overline\Omega_i$. We define cell membranes $\Gamma_i=\partial\Omega_i/\partial\Omega$ for $i=0,\ldots,N$ and the full membrane $\Gamma=\bigcup_{i=0}^{N}\Gamma_i$. Moreover, we define possible common interfaces, as $\Gamma_{ij}=\partial\Omega_i\cap\partial\Omega_j$, for $i,j=0,\ldots,N$ and $i\neq j$, satisfying  $\Gamma_{ij}\subseteq\Gamma_i$. We refer to Fig.~\ref{fig:EMI} for an example of this geometrical setting. We remark that in neuroscience applications a thin layer of extracellular media is always present between different cells, i.e. $\Gamma_{ij}=\varnothing$. On the other hand, in cardiac tissue modeling, different cells are often in direct contact, through so-called \textit{gap junction}, with a non-empty $\Gamma_{ij}$. We explore both applications, theoretically and experimentally. 

We consider the following stationary problem for extra- and intracellular potentials $u_i:\overline{\Omega}_i\to\mathbb{R}$, for $i=0,\ldots,N$:
\begin{align}
    & -\nabla \cdot (\sigma_i\nabla u_i) = 0 && \quad \mathrm{in} \quad\Omega_i,\label{eq::EMI_2}\\
    & -\sigma_i\nabla u_i\cdot\mathbf{n}_i=\sigma_j\nabla u_j\cdot\mathbf{n}_j=I_{ij} && \quad \mathrm{in} \quad \Gamma_{ij},\,\forall j\neq i, \label{eq::EMI_3} \\
    & v_{ij} = u_i-u_j  && \quad \mathrm{in} \quad \Gamma_{ij},\, \forall j\neq i, \label{eq::EMI_4} \\
    & \tau I_{ij} = v_{ij} + f_{ij} && \quad \mathrm{in} \quad \Gamma_{ij},\,\forall j\neq i, \label{eq::EMI_5}
\end{align}
where $\sigma_i\in\mathbb{R}_+$ is the conductivity in $\Omega_i$, $\mathbf{n}_i$ is the outer normal of $\partial\Omega_i$, $I_{ij}:\Gamma_{ij}\to\mathbb{R}$ is a membrane current, and $f_{ij}\in L^2(\Gamma_{ij})$ is a known membrane source satisfying $f_{ij}=-f_{ji}$, e.g. in the form 
\begin{equation}\label{eq::f_ij}
    f_{ij} = \left\{\begin{array}{lr}
        g, & \quad\text{for } i < j,\\
        -g, & \quad\text{for } i > j,        \end{array}\right.
\end{equation}
for a source function $g\in L^2(\Gamma)$.
The physical parameter $\tau=C_M^{-1}\Delta t\in\mathbb{R}_+$ includes the membrane capacitance $C_M$ and a discrete time step $\Delta t$. Essentially, the EMI problem consists of $N+1$ homogeneous Poisson problems coupled at the interfaces $\Gamma_i$, where the solution can be discontinuous, via a Robin-type condition \eqref{eq::EMI_3}-\eqref{eq::EMI_5}, 
In particular, eq. \eqref{eq::EMI_5} is obtained by an implicit-explicit (IMEX) time discretization of the membrane current, given a capacitive and an ionic contribution:
$$I_{ij} =I_{ij,\mathrm{cap}} + I_{ij,\mathrm{ion}} =C_M\frac{\partial v_{ij}}{\partial t} + I_{\mathrm{ion}}(v_{ij}).$$
The capacitive contribution $I_{ij,\mathrm{cap}}$ satisfies the capacitor current-voltage relation and $I_{\mathrm{ion}}:\mathbb{R}\to\mathbb{R}$ is an ionic current subject to further modelling, i.e. a possibly non-linear reaction term. Given the initial condition for the trans-membrane voltage $v_\mathrm{in}:\Gamma\to\mathbb{R}$, we have $\tau^{-1}(v_{ij}-v_\mathrm{in})=I_{ij}-I_{\mathrm{ion}}(v_\mathrm{in})$ for $0\leq i<j\leq N$ (considering $v_{ij}=-v_{ji}$ and $I_{ij}=-I_{ji}$), and define the source term in \eqref{eq::f_ij} as
\begin{equation}\label{eq::f}
f= \tau I_{\mathrm{ion}}(v_\mathrm{in}) - v_\mathrm{in}.
\end{equation}
 We close the EMI problem with homogeneous Neumann boundary conditions:
\begin{equation}\label{eq:bc}
     \nabla u_i \cdot \mathbf{n}_i = 0 \quad \mathrm{in} \quad \partial\Omega\cap\partial\Omega_i \quad \text{for}\quad i=0,\ldots,N.
\end{equation}
Because of the pure Neumann boundary condition, the solving potentials $u_i$ are defined up to an additive constant.
Uniqueness can be enforced via Lagrange multi\-pliers, pinning the solution in one point, or on the discrete level, providing a discrete nullspace to the solution strategy. 

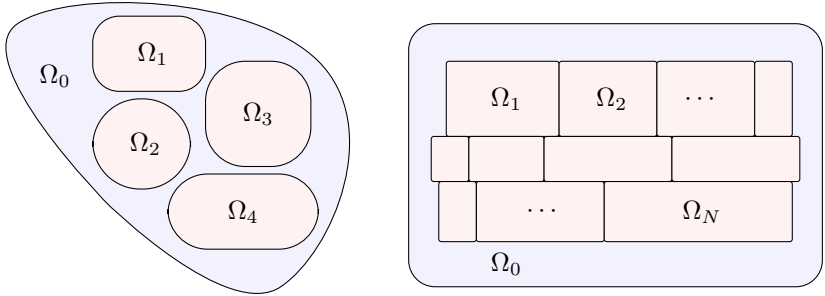
\begin{figure}
    \centering
    \begin{tikzpicture}

      \draw[fill=blue!5]  plot[very thick, smooth, tension=.9] coordinates {(-3.3,1.1) (-4.5,3.5) (-0.5,3) (-1,0) (-3.3,1.1)};

    \draw[rounded corners=10pt, fill=blue!5]
  (0.7,0) rectangle ++(5.5,3.5);

  \draw[rounded corners=18pt, fill=red!5]
  (-3.5,1.3) rectangle ++(1.3,1.2);

  \draw[rounded corners=15pt, fill=red!5]
  (-2,1.6) rectangle ++(1.4,1.4);

  \draw[rounded corners=15pt, fill=red!5]
  (-2.5,0.5) rectangle ++(2,1);

 \draw[rounded corners=10pt, fill=red!5]
  (-3.5,2.6) rectangle ++(1.5,1);

 \draw[rounded corners=1pt, fill=red!5]
  (1.2,2) rectangle ++(1.5,1);
 \draw[rounded corners=1pt, fill=red!5]
  (2.7,2) rectangle ++(1.3,1);
  \draw[rounded corners=1pt, fill=red!5]
  (4,2) rectangle ++(1.3,1);
  \draw[rounded corners=1pt, fill=red!5]
  (5.3,2) rectangle ++(0.5,1);
  
  \draw[rounded corners=1pt, fill=red!5]
  (1,1.4) rectangle ++(0.5,0.6);
  \draw[rounded corners=1pt, fill=red!5]
  (1.5,1.4) rectangle ++(1,0.6);
  \draw[rounded corners=1pt, fill=red!5]
  (2.5,1.4) rectangle ++(1.7,0.6);
  \draw[rounded corners=1pt, fill=red!5]
  (4.2,1.4) rectangle ++(1.7,0.6);

    \draw[rounded corners=1pt, fill=red!5]
  (1.1,0.6) rectangle ++(0.5,0.8);
    \draw[rounded corners=1pt, fill=red!5]
  (1.6,0.6) rectangle ++(1.7,0.8);
  \draw[rounded corners=1pt, fill=red!5]
  (3.3,0.6) rectangle ++(2.5,0.8);

    \node at (-4,2.8) {$\Omega_0$};
    \node at (-2.7,3.1) {$\Omega_1$};
    \node at (-2.8,1.9) {$\Omega_2$};
    \node at (-1.3,2.3) {$\Omega_3$};
    \node at (-1.5,1) {$\Omega_4$};

    \node at (2,0.3) {$\Omega_0$};
    \node at (2,2.5) {$\Omega_1$};
    \node at (3.4,2.5) {$\Omega_2$};
    \node at (4.6,2.5) {$\cdots$};
    \node at (2.5,1) {$\cdots$};
    \node at (4.6,1) {$\Omega_N$};

    \end{tikzpicture}
    \caption{Examples of geometrical setting for the nervous system (left) and cardiac tissue (right) for $d=2$.}
    \label{fig:EMI2}
\end{figure}

We investigate two cases of interest in the context of excitable tissues in physiologi\-cal applications, cf. Figure~\ref{fig:EMI2}, with final remarks on the general problem:

\begin{itemize}
    \item[(i)] A relevant setting to model nervous system cells. In this case, a layer of extracellular media is always present between different cells (e.g. neurons and glial cells, see for example \cite{sterratt2023principles}), i.e. $\Gamma_{i}=\Gamma_{i,0}$ for all $i=1,\ldots,N$ or, equivalently, $\Gamma_{ij}=\varnothing$ if $i,j\neq0$.
    \item[(ii)] A setting for cell-by-cell cardiac modeling. We consider the geometry of myocardiocytes, cardiac muscles cells which are connected by gap junctions or intercaleted discs, with a somehow regular pattern \cite{kleber2004basic}.
\end{itemize}

\section{Weak formulation and discrete operators}\label{sec:approx}

The EMI problem can be weakly formulated in various ways, depending on the unknowns of interest. We refer to \cite{EMI_book} for a broad discussion on various formulations (including the so-called mixed ones). As it could be expected from the structure of \eqref{eq::EMI_2}-\eqref{eq::EMI_5}, all formulations and corresponding discretizations give rise to block operators, with different blocks corresponding to $\Omega_i$, and $\Gamma_{ij}$.

We use a so-called \textit{single-dimensional} formulation and the corresponding discrete operators. In this setting, the weak form depends only on bulk quantities $u_i$  since the current term $I_{ij}$ is replaced by:
$$ I_{ij} = \tau^{-1}(u_i-u_j+f_{ij}),$$
according to equations~\eqref{eq::EMI_4}-\eqref{eq::EMI_5}.

After substituting the expression for $I_{ij}$ in \eqref{eq::EMI_3}, assuming the solution $u_i$ to be sufficiently regular over $\Omega_i$, we multiply the PDEs in \eqref{eq::EMI_2} by test functions $\phi_i\in V_i(\Omega_i)$, with $V_i$ a sufficiently regular Sobolev space with elements satisfying the boundary conditions in \eqref{eq:bc}; in practice $V_i=H^1(\Omega_i)$ is a standard choice. After integrating over $\Omega_i$ and applying integration by parts, using the normal flux definition \eqref{eq::EMI_3}, we have the following variational formulation for $i=0,\ldots,N$: 
\begin{align*}
& -\int_{\Omega_i} \nabla\cdot(\sigma_i\nabla u_i)  \phi_i\,\mathrm{d}\xx = 0,\\
    & \sigma_i\int_{\Omega_i} \nabla u_i \cdot\nabla \phi_i\,\mathrm{d}\xx - \int_{\partial \Omega_i}\sigma_i\nabla u_i\cdot \mathbf{n}_i\phi_i\,\mathrm{d}s  = 0, \\
     & \sigma_i\int_{\Omega_i} \nabla u_i \cdot\nabla \phi_i\,\mathrm{d}\xx + \sum_{j\neq i}\int_{\Gamma_{ij}}I_{ij}\phi_i\,\mathrm{d}s  = 0, \\
     & \tau\sigma_i\int_{\Omega_i} \nabla u_i \cdot\nabla \phi_i\,\mathrm{d}\xx + \sum_{j\neq i}\int_{\Gamma_{ij}}(u_i-u_j+f_{ij})\phi_i\,\mathrm{d}s = 0.
\end{align*}
Defining $\tau_i=\tau\sigma_i$, the weak EMI problem finally reads: for $i=0,\ldots,N$, find $u_i\in V_i(\Omega_i)$ such that
\begin{equation}
     \tau_i\int_{\Omega_i} \nabla u_i \cdot\nabla \phi_i\,\mathrm{d}\xx + \int_{\Gamma_{i}}u_i\phi_i\,\mathrm{d}s - \sum_{j\neq i}\int_{\Gamma_{ij}}u_j\phi_i\,\mathrm{d}s = -\sum_{j\neq i}\int_{\Gamma_{ij}}f_{ij}\phi_i\,\mathrm{d}s, \label{eq::EMI_weak1}
\end{equation}
for all test functions $\phi_i\in V_i(\Omega_i)$. We refer to \cite[Section 6.2.1]{EMI_book} for boundedness and coercivity results for this formulation.
According to \eqref{eq::f_ij}, the right hand side of \eqref{eq::EMI_weak1} can be rewritten as
\begin{equation}\label{eq::f_global}
    - \sum_{j\neq i}\int_{\Gamma_{ij}}f_{ij}\phi_i\,\mathrm{d}s = -\sum_{i<j}\int_{\Gamma_{ij}}g\phi_i\,\mathrm{d}s + \sum_{i>j}\int_{\Gamma_{ij}}g\phi_i\,\mathrm{d}s. 
\end{equation}

For each subdomain $\Omega_i$ we construct a conforming tessellation $\mathcal{T}_i$. We then introduce a yet unspecified discretization via finite element basis functions $\{\phi_{i,j}\}_{j=1}^{n_i}$ (e.g. Lagrangian elements of order $p\in\mathbb{N}$) for $V_{i,h}\subset V_i$:
$$V_{i,h}=\mathrm{span}\left(\{\phi_{i,k}\}_{k=1}^{n_i}\right), \quad \mathrm{and} \quad u_i(\xx)\approx u_{i,h}(\xx)=\sum_{k=1}^{n_i}u_{i,k}\phi_{i,k}(\xx),$$
with $n_i\in\mathbb{N}$ denoting the number of degrees of freedom in the corresponding sub\-domain $\Omega_i$ and $u_{i,k}\in\mathbb{R}$ the unknown coefficients.
Even if the presented theory is general, in the numerical experiments, we consider the natural choice of tessellations matching at each interface $\Gamma_i$.
%
From \eqref{eq::EMI_weak1} we define the following discrete operators for $i=0,\ldots,N$: bulk Laplacians
\begin{equation}    A_i=\left[\int_{\Omega_i}\nabla\phi_{i,\ell}(\xx)\cdot\nabla\phi_{i,k}(\xx)\,\mathrm{d}\xx\right]_{\ell,k=1}^{n_i}\in\mathbb{R}^{n_i\times n_i},
\end{equation}
and membrane mass matrices:\begin{equation}
    M_{i}=\left[\int_{\Gamma_i}\phi_{i,\ell}(\xx)\phi_{i,k}(\xx)\,\mathrm{d}s\right]_{\ell,k=1}^{n_i}\in\mathbb{R}^{n_i\times n_i}.
\end{equation}
We then define coupling matrices for $i,j=0,\ldots,N$ and $i\neq j$:
\begin{equation}
    B_{i,j}=-\left[\int_{\Gamma_{ij}}\phi_{i,\ell}(\xx)\phi_{j,k}(\xx)\,\mathrm{d}s\right]_{(\ell,k)=(1,1)}^{(n_i,n_j)}\in\mathbb{R}^{n_i\times n_j},
\end{equation}
and membrane source vectors
\begin{equation}
    \ff_{i}=-\left[\sum_{j \neq i}\int_{\Gamma_{ij}}f_{ij}(\xx)\phi_{i,k}(\xx)\,\mathrm{d}s\right]_{k=1}^{n_i}\in\mathbb{R}^{n_i}.
\end{equation}
Let us remark that the matrices $M_i$ and $B_{ij}$ are possibly both low rank, since have non-zero entries only corresponding to membranes $\Gamma_i$ degrees of freedom, cf. eq.~\eqref{eq:EMI_matrices}. We also note that $B_{i,j}=B^T_{j,i}$.
We write the symmetric linear system of size $n = \sum_{i={0}}^{N}n_i$ corresponding to \eqref{eq::EMI_weak1}:
\begin{equation}\label{eq::system}
	\begin{bmatrix}
 D_0 & B_{0,1} & B_{0,2} & \cdots & B_{0,N}  \\
	B_{1,0} &  D_1 & B_{1,2} & \ddots & \vdots  \\
    B_{2,0} &  B_{2,1} & D_2  & \ddots & B_{2,N}  \\
    \vdots &  \ddots & \ddots & \ddots & B_{N-1,N} \\    
    B_{N,0}  & \cdots & B_{N,N-2} & B_{N,N-1} & D_N  \\
	
	\end{bmatrix}
	\begin{bmatrix}
		\uu_0 \\
        \uu_1 \\
        \uu_2 \\
        \vdots \\
        \uu_{N-1} \\
        \uu_N \\		
	\end{bmatrix} =
	\begin{bmatrix}
		\ff_0 \\
        \ff_1 \\
        \ff_2 \\
        \vdots \\
        \ff_{N-1} \\
        \ff_N \\		
	\end{bmatrix}
 \quad \Longleftrightarrow \ \  \mathcal{A}_{\nn}\uu=\ff,
\end{equation}
with the diagonal blocks defined by
\begin{equation}\label{Di-expression}
D_i = \tau_i A_i + M_i,\qquad \text{for} \quad i=0,\ldots,N,
\end{equation}
and $\uu_i=[u_{i,0},u_{i,1},\ldots,u_{i,n_i}]\in\mathbb{R}^{n_i}$ collecting the finite element coefficients correspon\-ding to $\Omega_i$.
We label the coefficient matrix in \eqref{eq::system} as $\mathcal{A}_{\nn}$, {where $\nn$ is the vector $\nn = (n_0,n_1,\ldots,n_N)$ containing the sizes of the matrices $D_0,D_1,\ldots,D_N$. By summing the entries of $\nn$, we obtain the size \( n = \sum_{i=0}^N n_i \) of $\mathcal{A}_{\nn}$. Notice that if we fix the number of cells $N$, the domains, the type of tessellation, and the type of finite element basis functions, and only the mesh size varies, then the number \( n \) uniquely determines the vector $\nn$. In this setting, we can refer to the corresponding coefficient matrix sequence, as \(\{\mathcal{A}_{\nn}\}_n\).}
Furthermore, notice that each $n_i$, for $i=0,\ldots,N$, is also determined by various approximation and problem parameters as the degree $p$, the regularity $k$ of the global approximated solution (e.g. $k=0$ for standard Lagrangian finite elements, $k=-1$ for discontinuous Galerkin, and $k=p-1$ in the isogeometric analysis case) and on the dimensionality $d$ of the domain.
%
%
Moreover, making the entries corresponding to $\Gamma_i$ denoted explicitly (with subscript $\Gamma$), for $i=0,\ldots,N$, we can write
\begin{equation}\label{eq:EMI_matrices}
D_i=
	\begin{bmatrix}
		D_{i,\mathrm{b}} & D_{i,\mathrm{c}}  \\
		 D_{i,\mathrm{c}}^T & D_{i,\Gamma}
  \end{bmatrix},
  \qquad
  B_{i,j}=
  \begin{bmatrix}
		0 & 0  \\
		 0 & B_{ij,\Gamma}
  \end{bmatrix},
  \qquad
  M_i=
  \begin{bmatrix}
		0 & 0  \\
		 0 & M_{i,\Gamma}
  \end{bmatrix},
\end{equation}
where $D_{i,\mathrm{b}}$ are bulk blocks, $D_{i,\mathrm{c}}$ are coupling blocks, and $D_{i,\Gamma}, B_{ij,\Gamma}$ and $M_{i,\Gamma}$ are blocks corresponding to $\Gamma$, with size $n_{i, \Gamma}$. For notational convenience, we define the total number of intracellular and membrane points,
$$n_\mathrm{in}=\sum_{i=1}^Nn_i, \qquad n_\Gamma=\sum_{i=1}^Nn_{i,\Gamma}.$$
\subsection{Nervous system setting}\label{sec:nervous_system}

In this setting, we have $\Gamma_{ij}=\varnothing$ for $i,j\neq0$, i.e. the extra-diagonal blocks vanish for $i>0$, i.e. $B_{i,j}=0$ for $i,j>0$ and, with abuse of notation, the coefficients matrix $\mathcal{A}_{\nn}$ has the following block arrowhead form:
\begin{equation}\label{eq:sys_brain}
\mathcal{A}_{\nn}= 
\begin{bmatrix}
 D_0 & B_1 & B_2 & \cdots & B_N  \\
	B_1^T &  D_1 & 0 & \cdots & 0  \\
    B_2^T &  0 & D_2  & \ddots & \vdots  \\
    \vdots &  \vdots & \ddots & \ddots & 0 \\    
    B_N^T  & 0 & \cdots & 0 & D_N 	
	\end{bmatrix},
\end{equation}
where $B_j\equiv B_{0,j}=B_{j,0}^T$ for notational convenience. 
Analytical inverses of block arrowhead matrices \cite{stanimirovic2019inversion}, depending only on sub-blocks inversions, can be obtained via the Sherman–Morrison–Woodbury (SMW) formula (a.k.a. Woodbury matrix identity), depending on a decomposition in the form:
\begin{equation}\label{matr:a-nn}
\mathcal{A}_{\nn}=D_{\nn}+U_{\nn}V_{\nn},
\end{equation}
with $D_{\nn},U_{\nn},V_{\nn}$ conformable matrices. According to \eqref{eq::system}, in the EMI case the following decomposition is a natural choice:
\begin{equation}\label{matr:d-nn etc}
 D_{\nn}=
 \begin{bmatrix}
 D_0 & & & &   \\
&  D_1 &  & &   \\
    & &  \ddots &  & &   \\
    & & &  D_N &     \\
	\end{bmatrix}, \,
 U_{\nn} =
 \begin{bmatrix}
 I_{n_0} & 0 \\
 0 & B_1^T \\
 \vdots & \vdots \\
 0 & B_N^T \\
	\end{bmatrix},\,
 V_{\nn} =
 \begin{bmatrix}
 0 & B_1 & \cdots & B_N \\
 I_{n_0} & 0 & \cdots & 0 \\
	\end{bmatrix}.
\end{equation}
Crucially, the $D_1,\ldots,D_N$ blocks are singular, since they correspond to sub-problems with unspecified boundary conditions (i.e. a rank-one deficiency), and therefore the block diagonal matrix $D_{\nn}$ is also singular, with $\mathrm{rank}(D_{\nn}) = n-N$, while $U_nV_n$ has rank $2n_\Gamma$ (by numerical inspection). Thus, to apply the SMW formula, given $\epsilon>0$, we define a full rank correction of $D_{\nn}$ as 
$$D_{\nn,\epsilon} = D_{\nn} + \epsilon I_n,$$
 and the corresponding block arrowhead matrix
$$\mathcal{A}_{\nn,\epsilon} =  D_{\nn,\epsilon} + U_{\nn}V_{\nn},$$
with $\lim_{\epsilon\to 0}\mathcal{A}_{\nn,\epsilon} = \mathcal{A}_{\nn}$.
Notice that the term $\epsilon I_n$ could also be replaced by any term $\epsilon E_n$ such that $D_{\nn} + \epsilon E_n$ is invertible for any $\epsilon>0$.
We can then apply the SMW formula to obtain $\mathcal{A}_{\nn,\epsilon}^{-1}$ as
\begin{equation}\label{eq::SMW-eps}
    \mathcal{A}_{\nn,\epsilon}^{-1} = D_{\nn,\epsilon}^{-1} - D_{\nn,\epsilon}^{-1}V_{\nn}(I_{2n_0}+ V_{\nn}D_{\nn,\epsilon}^{-1}U_{\nn})^{-1}U_{\nn}D_{\nn,\epsilon}^{-1},
\end{equation}
which depends on $D_{\nn,\epsilon}^{-1}$, i.e. $N$ independent inversions $D_0^{-1},\ldots,D_N^{-1}$, and the solution of a linear system of size $2n_0$. Moreover, given $\lim_{\epsilon\to 0}\mathcal{A}_{\nn,\epsilon}^{-1} = \mathcal{A}_{\nn}^{-1}$, the form of $\mathcal{A}_{\nn,\epsilon}^{-1}$ can be used to design preconditioners, cf. equation~\eqref{eq::P}.
Alternatively, we introduce the following decomposition of $\mathcal{A}_{\nn}$ as 
$$\mathcal{A}_{\nn} =  D_{\nn} + E_{\nn} - E_{\nn} + U_{\nn}V_{\nn} ,$$
 with the rank $N$ correction:
\begin{equation*}
E_{\nn} =
\begin{bmatrix}
 0 & 0 & 0 & 0 \\
 0 & \mathbf{e}^T_{1,1}\mathbf{e}_{1,1} & 0 & 0 \\
 0 & 0 & \ddots & 0 \\
 0 & 0 & 0 & \mathbf{e}^T_{1,N}\mathbf{e}_{1,N}
\end{bmatrix},
\end{equation*}
{where $\mathbf{e}_{1,i}=[1,0,\ldots,0]\in\mathbb R^{n_i}$}. We can then write
$$ \mathcal{A}_{\nn} = \widetilde{D}_{\nn} + \widetilde{U}_{\nn}\widetilde{V}_{\nn},$$
with $\widetilde{D}_{\nn} = D_{\nn} + E_{\nn}$, 
\begin{equation*}
    \widetilde{U}_{\nn} =
    \left[
    \begin{array}{c|ccccc}
 & 0 & 0 & 0 & 0 \\
 & \mathbf{e}^T_{1,1} & 0 & 0 & 0 \\
U_{\nn} & 0 & \mathbf{e}^T_{1,2} & 0 & 0  \\
 & 0 & 0 & \ddots & 0 \\
 & 0 & 0 & 0 & \mathbf{e}^T_{1,N} \\
\end{array}
    \right], \quad
    \widetilde{V}_{\nn} =
\left[
\begin{array}{ccccc}
 & & V_{\nn} & & \\
 \hline
0 & \mathbf{e}_{1,1} & 0 & 0 & 0 \\
0 & 0 & \mathbf{e}_{1,2} & 0 & 0  \\
0 & 0 & 0 & \ddots & 0 \\
0 & 0 & 0 & 0 & \mathbf{e}_{1,N} \\
\end{array}
\right],
\end{equation*}
where $\widetilde{D}_{\nn}$ has full rank and $\widetilde{U}_{\nn}\widetilde{V}_{\nn}$ is a $2n_\Gamma + N$ rank correction. Again, applying the SMW formula we obtain a closed equation for $\mathcal{A}_{\nn}^{-1}$:
\begin{equation}\label{eq::SMW}
    \mathcal{A}_{\nn}^{-1} = \widetilde{D}_{\nn}^{-1} - \widetilde{D}_{\nn}^{-1}\widetilde{V}_{\nn}(I_{2n_0+N}+ \widetilde{V}_{\nn}\widetilde{D}_{\nn}^{-1}\widetilde{U}_{\nn})^{-1}\widetilde{U}_{\nn}\widetilde{D}_{\nn}^{-1}.
\end{equation}
We emphasize that both formulae (\ref{eq::SMW-eps}) and (\ref{eq::SMW}) are of concrete algorithmic interest when either $2n_0$ or $2n_0+N$ are negligible compared with global matrix-size $n$. However, the size $n_0$, i.e. the number of extra-cellular degrees of freedom,  still corresponds to a $d$-dimensional discrete Laplacian block ($D_0$) with $d=2,3$. Hence, unless one uses a symbol approach as in the Schur complement computation in \cite{NS1-elonged} by means of circulant or $\tau$ based approximations, these formulae do not offer a clear efficient solution strategy. 


\subsection{Cardiac setting}\label{sec:cardiac_setting}
In the cardiac setting, cells can be in direct contact, i.e. with a non-empty $\Gamma_{ij}$, and the block  structure in \eqref{eq:sys_brain} is lost. 
In particular, the structure of $\mathcal{A}_{\nn}$ remains general but block sparse. More precisely, the resulting sparsity pattern has a multilevel band shape: for instance, a two-level tridiagonal with tridiagonal blocks is encountered when considering a structured cell reticulum for $d=2$. In that specific setting, the nondiagonal blocks are essential for determining the global spectral distribution and this is the reason why a pure block diagonal preconditioner might not be effective in this case, cf. Section~\ref{sec:num}.

\section{Spectral analysis}\label{sec:spectral}

In this section, we study the spectral distribution of the matrix sequence $\{\mathcal{A}_{\nn}\}_{n}$ in \eqref{eq:sys_brain} under various assumptions for determining the global behaviour of the eigenvalues of $\mathcal{A}_{\nn}$ as the matrix-size $n = \sum_{i=0}^{N}n_i$ tends to infinity
and for a fixed number of cells $N$, potentially large: the limit case of $N\to\infty$ is briefly discussed in Remark \ref{rem:N infinity}. The spectral distribution is given by a smooth function called the (spectral) symbol as it is customary in the Toeplitz and Generalized Locally Toeplitz (GLT) setting  \cite{GS-I,GS-II,block-glt-1D,block-glt-dD}.

We consider in detail the setting describing cells in the nervous system, with $\mathcal{A}_{\nn}$ given by \eqref{eq:sys_brain}. A short discussion on the cardiac setting is contained in Remark \ref{rem:cardiac-trid}.
Regarding the needed notions,  we give the formal definition of Toeplitz structures, eigenvalue (spectral) and singular value distribution, few connected basic tools, and finally, we provide the specific analysis of EMI matrix sequences under a variety of assumptions.

\subsection{Toeplitz structures, spectral symbol, GLT tools}\label{ssec:GLTbackground}

We provide the definition of block Toeplitz sequences associated with a matrix-valued Lebesgue inte\-grable function over $[-\pi,\pi]^d$, $d\ge 1$. Subsequently, we introduce the notion of eigenvalue (spectral) and singular value distribution, and we report few tools taken from the 
relevant literature.

\begin{definition}\label{def_not3}{\rm [Toeplitz sequences (generating function of)]}
Denote by $f$ a $d$-variate complex-valued integrable function,
defined over the domain $Q^d=[-\pi,\pi]^d, d\ge 1$, with $d$-dimensional Lebesgue measure $\mu_d (Q^d)=(2\pi)^d$. Denote by $f_k$ the Fourier coefficients of $f$,
$$
f_k = \frac{1}{(2\pi)^d}\int\limits_{Q^d}f(\theta)e^{-i\,(k,\theta)}\,d\theta,\;
k=(k_1,\cdots,k_d)\in \ZZ^d, \; i^2=-1,
$$
where $\theta=(\theta_1,\cdots,\theta_d)$, $(k,{ \theta})=\sum_{j=1}^{d}k_j \theta_j$. {By following the multi-index notation in \cite{tyrt}[Section 6], with each
$f$ we can associate a sequence of Toeplitz matrices $\{T_\nu\}_\nu$, with $\nu=(\nu_1,\cdots,\nu_d)$ being a multi-index, with the formula
$$
T_\nu = \{ f_{k-\ell}\}_{k,\ell=\mathbf{e}^T}^{\nu} \in \mathbb{C}^{N(\nu)\times N(\nu)},
$$
with $\mathbf{e}=[1,1,\cdots,1]\in \NN^d$ and $N(\nu)={\nu_1}\times\cdots\times \nu_d$.}

For $d=1$
$$
T_\nu =
\begin{bmatrix}
f_0&f_{-1}&\cdots&f_{-\nu+1}\\
f_{1}&f_0&\ddots&\vdots\\
\vdots&\ddots&f_0&f_{-1}\\
f_{\nu-1}&\cdots&f_1&f_0
\end{bmatrix},
$$
or for $d = 2$, i.e. the two-level case, and for example $\nu=(2,3)$, we have
$$
T_\nu =
\begin{bmatrix}
F_0&F_{-1}\\
F_{1}&F_0
\end{bmatrix},
\quad
F_k =
\begin{bmatrix}
f_{k,0}&f_{k,-1}&f_{k,-2}\\
f_{k,1}&f_{k,0}&f_{k,-1}\\
f_{k,2}&f_{k,1}&f_{k,0}
\end{bmatrix},\ \ \ k=0,\pm 1.
 $$
The function $f$ is referred to as the \textit{generating function} (or the
\textit{symbol} of) $T_\nu$.  Using a more compact notation, we say that the
function $f$ is the generating function of the whole sequence $\{T_\nu\}_\nu$ and we
write $T_\nu = T_\nu(f)$.

If $f$ is $d$-variate, $\mathbb{C}^{s_1\times s_2}$ matrix-valued, and integrable over  $Q^d$, $d,s_1,s_2\ge 1$, i.e.
$f\in L^1(Q^d,s_1\times s_2)$, then we can define the
Fourier coefficients of $f$ in the same way (now $f_k$ is a matrix of size
$s_1\times s_2$) and consequently  $T_\nu = \{ f_{k-\ell}\}_{k,\ell=\mathbf{e}^T}^{\nu}
\in \mathbb{C}^{s_1N(\nu)\times s_2N(\nu)}$, then $T_\nu$ is a $d$-level block
Toeplitz matrix according to Definition 4.1 in \cite{benedusi2023EMI}.
If $s_1=s_2=s$ then we write $f\in L^1(Q^d,s)$.

As in the scalar case, the function $f$ is referred to as the \textit{generating function} of $T_\nu$. We say that the function $f$ is the generating
function of the whole sequence $\{T_\nu\}_\nu$, and we use the notation $T_\nu = T_\nu(f)$.
\end{definition}

\begin{definition}\label{def-distribution}
	Let $f:D\to\mathbb{C}^{s\times s}$ be a measurable matrix-valued function with eigenvalues $\lambda_i(f)$ and singular values $\sigma_i(f)$, $i=1,\ldots,s$. Assume that $D\subset \mathbb{R}^d$ is Lebesgue measurable with positive and finite Lebesgue measure $\mu_d(D)$. Assume that $\{A_n\}_n$ is a sequence of matrices such that ${\rm dim}(A_n)=d_n\rightarrow\infty$, as $n\rightarrow\infty$ and with eigenvalues $\lambda_j(A_n)$ and singular values $\sigma_j(A_n)$, $j=1,\ldots,d_n$.
	\begin{itemize}
		\item We say that $\{A_n\}_{n}$ is {\em distributed as $f$ over $D$ in the sense of the eigenvalues,} and we write $\{A_n\}_{n}\sim_\lambda(f,D),$ if
		\begin{equation}\label{distribution:eig}
		\lim_{n\to\infty}\frac{1}{d_n}\sum_{j=1}^{d_n}F(\lambda_j(A_n))=
		\frac1{\mu_d(D)} \int_D \frac1{s}\sum_{i=1}^sF(\lambda_i(f(t))) \, \mathrm{d}t,
		\end{equation}
		for every continuous function $F$ with compact support. In this case, we say that $f$ is the \emph{spectral symbol} of $\{A_{n}\}_{n}$.
		
		\item We say that $\{A_n\}_{n}$ is {\em distributed as $f$ over $D$ in the sense of the singular values,} and we write $\{A_n\}_{n}\sim_\sigma(f,D)$, if
		\begin{equation}\label{distribution:sv}
		\lim_{n\to\infty}\frac{1}{d_n}\sum_{j=1}^{d_n}F(\sigma_j(A_n))=
		\frac1{\mu_d(D)} \int_D \frac1{s}\sum_{i=1}^sF(\sigma_i(f(t))) \, \mathrm{d}t,
		\end{equation}
		for every continuous function $F$ with compact support. In this case, we say that $f$ is the \emph{singular value symbol} of $\{A_{n}\}_{n}$.
		\item The notion $\{A_n\}_{n}\sim_\sigma(f,D)$ applies also in the rectangular case where $f$ is $\mathbb{C}^{s_1\times s_2}$ matrix-valued. In such a case the parameter $s$ in formula (\ref{distribution:sv}) has to be replaced by the minimum between $s_1$ and $s_2$: furthermore $A_n\in \mathbb{C}^{d_n^{(1)}\times d_n^{(2)}}$ with $d_n$ in formula (\ref{distribution:sv}) being the minimum between $d_n^{(1)}$ and $d_n^{(2)}$. Of course the notion of eigenvalue distribution does not apply in a rectangular setting.
			\end{itemize}
\end{definition}
Throughout the paper, when the domain can be easily inferred from the context, we replace the notation $\{A_n\}_n\sim_{\lambda,\sigma}(f,D)$ with $\{A_n\}_n\sim_{\lambda,\sigma} f$.

\begin{remark}\label{rem:approx}
	If $f$ is smooth enough, an informal interpretation of the limit relation
	\eqref{distribution:eig} (resp. \eqref{distribution:sv})
	is that when $n$ is sufficiently large, the eigenvalues (resp. singular values) of $A_{n}$ can be subdivided into $s$ different subsets of the same cardinality. Then $d_n/s$ eigenvalues (resp. singular values) of $A_{n}$ can
	be approximated by a sampling of $\lambda_1(f)$ (resp. $\sigma_1(f)$)
	on a uniform equispaced grid of the domain $D$, and so on until the
	last $d_n/s$ eigenvalues (resp. singular values), which can be approximated by an equispaced sampling
	of $\lambda_s(f)$ (resp. $\sigma_s(f)$) in the domain $D$.
\end{remark}
\begin{remark}
    We say that $\{A_n\}_n$ is \textit{zero-distributed} in the sense of the eigenvalues if $\{A_n\}_n\sim_\lambda 0$. Of course, if the eigenvalues of $A_n$ tend all to zero for $n\to\infty$, then this is sufficient to claim that $\{A_n\}_n\sim_\lambda 0$.

\end{remark}
For Toeplitz matrix sequences, the following theorem due to Tilli holds, which genera\-lizes previous research along the last 100 years by Szeg\H{o}, Widom, Avram, Parter, Tyrtyshni\-kov, and Zamarashkin (see \cite{GS-II,block-glt-dD} and references therein).

\begin{theorem}{\rm \cite{TilliNota}}\label{szego-herm}
	Let $f\in L^1(Q^d,s_1\times s_2)$, then $\{T_{\nu}(f)\}_{{\nu}}\sim_\sigma(f,Q^d).$ If $s_1=s_2=s$ and if $f$ is a Hermitian matrix-valued function, then $\{T_{\nu}(f)\}_{{\nu}}\sim_\lambda(f,Q^d)$.
\end{theorem}

The following theorem is useful for computing the spectral distribution of a sequence of Hermitian matrices. For the related proof, see \cite[Theorem 4.3]{curl-div} and \cite[Theorem 8]{NS1-elonged}. Here, the conjugate transpose of the matrix $X$ is denoted by $X^*$.

\begin{theorem}{\rm \cite[Theorem 4.3]{curl-div}}\label{th:extradimensional}
	Let $\{A_n\}_n$ be a sequence of matrices, with $A_n$ Hermitian of size $d_n$, and let $\{P_n\}_n$ be a sequence such that $P_n\in\mathbb C^{d_n\times\delta_n}$, $P_n^*P_n=I_{\delta_n}$, $\delta_n\le d_n$ and $\delta_n/d_n\to1$ as $n\to\infty$. Then $\{A_n\}_n\sim_{\lambda}f$ if and only if $\{P_n^*A_nP_n\}_n\sim_{\lambda}f$.
\end{theorem}
With the notations of the result above, the matrix sequence $\{P_n^*A_nP_n\}_n$ is called a compression of $\{A_n\}_n$ and the single matrix $P_n^*A_nP_n$ is called a compression of $A_n$.

In what follows we take into account a crucial fact that is often neglected: the generating function of a Toeplitz matrix sequence and even more the spectral symbol of a given matrix sequence is not unique, except for the trivial case of either a constant generating function or a constant spectral symbol. In fact, here we report and generalize \cite[Remark 1.3]{Schur-CMAME} and the discussion below \cite[Theorem 3]{Symbols-SE}.

\begin{remark}
\noindent We remark that the presented tools are general and can be applied to matrix sequences stemming from a variety of discretization schemes such as isogeome\-tric analysis or finite volumes, in the spirit of the sections of books \cite{GS-II,block-glt-dD} dedicated to applications and of the exposition paper \cite{garoni2019}.
\end{remark}


\subsection{Symbol analysis}\label{ssec:symbol}

Here, we state and prove three results which hold under reasonable assumptions. The first is provided in the maximal generality, while the second and the third can be viewed as special cases of the first.
The results extend substantially those in \cite{benedusi2023EMI}[Section 4.2] where a single cell in the EMI model is considered, while the other level of generalization concerns the larger range of relations among all the various scaled time-steps and space-steps, $\tau_i=\sigma_i\tau, h_i$, $i=0,\ldots,N$.

\begin{theorem}\label{th:general}
For $i=0,\ldots,N$ assume that
\begin{align*}
& n_{i,\Gamma}  =  o(\min\{n_i,n_0\})  \quad \mathrm{for} \quad n_{i,\Gamma},n_i\to\infty, \\
& \lim_{n_i\rightarrow\infty }  \frac{n_i}{n} = r_i \in (0,1), \\  
& n=n_0+n_\mathrm{in}=\sum_{i=0}^Nn_i,\ \  \sum_{i=0}^N r_i=1.
\end{align*}
Also assume
\begin{equation}\label{distrAi}
    \{D_i\}_n\sim_\lambda (f^i,D^i),
\end{equation}
with $D^i\subset \mathbb{R}^{k_i}$, given $f^i$ the $i$th cell symbol.
Taking into consideration the matrix structures given in (\ref{matr:a-nn})-(\ref{matr:d-nn etc}), it follows that
\begin{eqnarray*}
\{D_{\nn}\}_n & \sim_\lambda & (g,[0,1]\times D),\\
\{\mathcal{A}_{\nn}\}_n & \sim_\lambda & (g,[0,1]\times D), \\
\{\mathcal{A}_{\nn}-D_{\nn}\}_n & \sim_\lambda & (0,[0,1]\times D),
\end{eqnarray*}
where $D=D^0\times D^1\times D^2 \times \cdots \times D^N$ with
\[
g(x,t_0,t_1,\ldots,t_N)= \sum_{i=0}^N f^i(t_i)\psi_{[\hat{r}_{i-1} ,\hat{r}_{i}]}(x),
\]
$x\in [0,1]$, $t_i\in D^i$ and $\hat{r}_{-1}=0$, $\hat{r}_{i}=\hat{r}_{i-1}+r_i$, $\hat{r}_{N}=1$,
$\psi_Z$ denoting the characteristic function of the set $Z$. 
\end{theorem}
\ \\
{\bf Proof}\,
As a preliminary observation, the employed Galerkin approach, independently of the specific method, and the structure of the equations imply that all the involved matrices are real and symmetric. Hence the symbols $f^i$, $i=0,\ldots,N$, are necessarily Hermitian matrix-valued. For the sake of notational simplicity we assume that all $f^i$ take values into $\mathbb{C}^{s\times s}$ with a fixed $s$ independent of the various parameters: the latter can be forced without loss of generality thanks to Remark 4.3 and Remark 4.6 in \cite{benedusi2023EMI}, where the case of different matrix-sizes in the symbols can be reduced to the same matrix-size with a trick based on the non-uniqueness of the spectral symbols.
Taking into account Definition \ref{def-distribution} and the assumptions, we have
\begin{eqnarray} 
\lim_{n_i\to\infty}\frac{1}{n_i}\sum_{j=1}^{n_i}F(\lambda_j(D_i)) & = &
		\frac1{\mu_{k_i}(D^i)} \int_{D^i} \frac1{s}\sum_{m=1}^sF(\lambda_m(f^i(t_i))) \,\mathrm{d}t_i,		\label{distr-i}
\end{eqnarray}
for $i=0,\ldots,N$ and for any continuous function $F$ with bounded support. Now we consider the $(N+1)\times (N+1)$ block diagonal matrix as in
(\ref{matr:d-nn etc})
\[
D_{\nn}={\rm diag}_{i=0}^N(D_i)
\]
and a generic continuous function $F$ with bounded support. By defining 
\[
\Delta(D_{\nn},F) =   \frac{1}{n}\sum_{i=0}^N\sum_{j=1}^{n_i}F(\lambda_j(D_i)),
\]
by exploiting the block diagonal structure of $D_{\nn}$, we infer 

\begin{eqnarray*}
\Delta(D_{\nn},F) & = &   \sum_{i=0}^N\frac{n_i}{n} \frac{1}{n_i}\sum_{j=1}^{n_i}F(\lambda_j(D_i)).
\end{eqnarray*}

As a consequence, according to \eqref{distr-i} and using the spectral symbol notion, we obtain
that the limit of $\Delta(D_{\nn},F)$ for $n_i\to \infty$ with $i=0,\ldots,N$ exists and 
\begin{equation}\label{eq: ci-siamo}
{\lim_{n_i\to\infty} \Delta(D_{\nn},F)}=\sum_{i=0}^N \frac{r_i}{\mu_{k_i}(D^i)} \int_{D^i} \frac1{s}\sum_{m=1}^sF(\lambda_m(f^i(t_i))) \,\mathrm{d}t_i.
\end{equation}
When looking at Definition \ref{def-distribution}, the quantity (\ref{eq: ci-siamo}) is not in the form of the righthand side of (\ref{distribution:eig}): in fact, we observe a summation with $N+1$ different terms. The difficulty is not serious and can be overcome by enlarging the space with a fictitious domain $[0,1]$ and interpreting the sum in (\ref{eq: ci-siamo}) as the global integral involving a step function. By leveraging the previous statement, we rewrite (\ref{eq: ci-siamo}) as
\begin{eqnarray*}
\frac{1}{\mu_{k+1}(\hat D)}
\int_0^1 \,\mathrm{d}x \int_{D} \frac1{s}\sum_{m=1}^sF {\left( 
\sum_{i=0}^N\lambda_m(f^i(t_i) \psi_{[\hat{r}_{i-1} ,\hat{r}_{i}]}(x))\right)
\mathrm{d}t_0\mathrm{d}t_1\ldots\mathrm{d}t_N},
\end{eqnarray*}
with $\hat D=[0,1]\times D$ and 
$k=\sum_{i=0}^N k_i$. As a consequence,
\[
\{D_{\nn}\}_n\sim_\lambda (g,[0,1]\times D)
\]
is complete.

For the study of $\mathcal{A}_{\nn}$ we notice that rank$(\mathcal{A}_{\nn}-D_{\nn})$ depends on the matrix-sizes associated with the {membrane} terms.
Therefore, by using the first assumption, rank$(\mathcal{A}_{\nn}-D_{\nn})\le 2\sum_{i=1}^N n_{i,\Gamma}  =  o(n)$, $n=\sum_{i=0}^N n_i$.
The latter is sufficient, by explicit computation, to claim that the related matrix-sequence is zero-distributed in the eigenvalue sense i.e.
$\{\mathcal{A}_{\nn}-D_{\nn}\}_n\sim_\lambda 0$. Since $\mathcal{A}_{\nn}= D_{\nn} + R_{\nn}$, $R_{\nn}=\mathcal{A}_{\nn}-D_{\nn}$, the statement in \cite[Exercise 5.3]{GS-I} implies directly that $\{\mathcal{A}_{\nn}\}_n$ and $\{D_{\nn}\}_n$ share the same distribution i.e.
$\{\mathcal{A}_{\nn}\}_n\sim_\lambda (g,[0,1]\times D)$ and the proof is concluded.
\ \hfill $\bullet$

The following two corollaries simplify the statement of Theorem \ref{th:general}, under special assumptions which are satisfied for few basic discretization schemes and when dealing with elementary domains. The argument behind the possibility of writing several (indeed infinitely many) spectral symbols relies on the non-uniqueness of the spectral symbol and on the rearrangement theory: for a discussion on the matter refer to \cite{Symbols-SE,barbarino2022constructive}.

\begin{corollary}\label{th:specific}
With the same notations and assumptions as in Theorem \ref{th:general}, assuming that $r_i= \frac{1}{N+1}$, $i=0,\ldots,N$ and $D^i=\tilde D\subset\mathbb{R}^{k}$, $i=0,\ldots,N$, we deduce that
\begin{eqnarray*}
\{D_{\nn}\}_n & \sim_\lambda & \left({\rm diag}(f^0,f^1,\ldots,f^N),\tilde D\right), \\
\{\mathcal{A}_{\nn}\}_n & \sim_\lambda &  \left({\rm diag}(f^0,f^1,\ldots,f^N),\tilde D\right),
\end{eqnarray*}
and $\{\mathcal{A}_{\nn}-D_{\nn}\}_n\sim_\lambda 0$.
\end{corollary}
{\bf Proof}\,
Under the given assumptions, the thesis of Theorem \ref{th:general} is equivalent to the main statement, 
taking into account Definition \ref{def-distribution} in the matrix-valued setting, i.e. relation (\ref{distribution:eig}) with diagonal matrix-valued spectral symbol and $s=N+1$.
For the last statement, $\{\mathcal{A}_{\nn}-D_{\nn}\}_n\sim_\lambda 0$, it is sufficient to follow verbatim the same reasoning as in Theorem \ref{th:general}.
\ \hfill $\bullet$

\begin{corollary}\label{th:specific-bis}
With the same notations and assumptions as in Theorem \ref{th:general}, assuming $f^i=f$, $D^i=\tilde D\subset\mathbb{R}^{k}$, $i=0,\ldots,N$, we deduce that
\[
\{D_{\nn}\}_n,\ \ \{\mathcal{A}_{\nn}\}_n\sim_\lambda  (f,\tilde D)
\]
and $\{\mathcal{A}_{\nn}-D_{\nn}\}_n\sim_\lambda 0$.

\end{corollary}
{\bf Proof}\,
The proof of the relation $\{\tilde A_n\}_n,\ \ \{A_n\}_n\sim_\lambda (f,\tilde D)$ follows directly from the limit displayed in (\ref{eq: ci-siamo}),
after replacing $f^i$ with $f$ and necessarily $D^i$ with $D$, $i=0,\ldots,N$. The rest is again obtained verbatim as in Theorem \ref{th:general}.
\ \hfill $\bullet$

Few remarks are in order.               
\begin{remark}\label{rem:N infinity}
According to the study in \cite{Symbols-SE,barbarino2022constructive}, the spectral symbol in the Weyl sense stated Definition \ref{def-distribution} is far from unique and in fact any rearrangement is still a symbol. Therefore we have infinitely many choices. For a fixed value of $N$, the expressions in Theorem \ref{th:general} and Corollary \ref{th:specific} are satisfactory, but they cannot be used as $N$ tends to infinity since either the domain 
$D^0\times D^1\times D^2 \times \cdots \times D^N$ changes dramatically with $N$ or the size of the symbol explodes with $N$. However, when considering the univariate nondecreasing rearrangement $\phi_N$ of the symbol and under mild assumptions, the limit as $N$ tends to infinity can be computed and this allows to understand the limit spectral symbol when the number of cells is extremely large, so allowing to treat very realistic situations.
\end{remark}
\begin{remark}\label{rem:cardiac-trid}
In the cardiac setting, due to the direct contact of the cells, the block  structure in \eqref{eq:sys_brain} is lost and the  important part of the matrix has a sparsity pattern of $d$-level tridiagonal type with $d\ge 2$. Hence the nondiagonal blocks are essential for determining the global spectral distribution and this is the reason why a pure block diagonal preconditioner is not robust with respect to the main parameters, since the related matrix-sequence does not capture the symbol of the discretization matrix-sequence. The analysis is not carried out in detail here and it will be the subject of further investigations.
\end{remark}

\subsection{Varying and adapting the symbol analysis}\label{ssec:symbol-bis}

The symbol analysis given in the previous section is now critically discussed in terms of the meaning of the assumptions. First, we should consider the significance of assuming the existence of the distribution symbols $f^i$, $D^i\subset \mathbb{R}^{k_i}$, $i=0,\ldots,N$, related to the various sequences $\{D_i\}_n$, associated with the blocks of $\{D_{\nn}\}_n$, namely, equation \eqref{distrAi}.
Looking at the structure of $D_i$, as reported in \eqref{Di-expression}, assumption \eqref{distrAi} holds if and only if
\[
    \lim_{h_i,\tau_i\to 0} \frac{\tau_i}{h_i^2} = c_i \in [0,\infty), \ \ \ \ i=0,\ldots,N.
\]
Indeed, the GLT analysis of stiffness matrices arising from various discretizations (see e.g. \cite{DG,benedusi2018space} for discontinuous Galerkin, \cite{GSS-fem,Schur-CMAME} for finite elements of any order, \cite{IgA-colloc,IgA-galerkin,garoni2019} for the isogeometric analysis both in Galerkin and collocation form) guarantees that there exist $f^i$ and $D^i$ such that
\begin{equation}\label{hp-distrib}
\left\{h_i^2 A_i\right\}_n \sim_\lambda (f^i, D^i), \ \ \ i=0,\ldots,N,
\end{equation}
and so we have
\[
    \left\{{\tau_i} A_i\right\}_n =
    \left\{\frac{\tau_i}{h_i^2} h_i^2 A_i\right\}_n \sim_\lambda (c_i f^i, D^i).
\]
Moreover, the rank of \(M_i\) in (\ref{Di-expression}) is equal to the rank of \(M_{i,\Gamma}\), which, in turn, is bounded by $n_{i,\Gamma}  =  o(\min\{n_i,n_0\})$ for $n_{i,\Gamma},n_i\to\infty$. Therefore, the related matrix sequence is zero-distributed in the eigenvalue sense, i.e., \(\{M_i\}_n \sim_\lambda 0\). Since \(D_i = \tau_i A_i + M_i\), applying the statement in \cite[Exercise 5.3]{GS-I} we get \(\{D_i\}_n \sim_\lambda (c_i f^i, D^i)\).

However, we note that if \( c_i = 0 \) for some \( i \), then the corresponding matrix sequences $\{D_i\}_n$ are zero-distributed in the eigenvalue sense. According to Theorem \ref{th:general}, this leads to the problem being severely ill-posed. Consequently, we choose to exclude this case from our practical consideration.

In the case where
\[
    \lim_{h_i,\tau_i\to 0} \frac{\tau_i}{h_i^2} = \infty, \ \ \ \ i=0,\ldots,N,
\]
or where the limit above holds for some of the values $i=0,\ldots,N$, we need to scale the matrices $D_i$ by a factor of $h_i^2/\tau_i$ for equation \eqref{distrAi} to hold. So, we perform the following scaling to the global matrix $\mathcal{A}_{\nn}$
\[
    \mathcal{A}_{\nn}^{\rm{s}} =
    {\rm diag}_{i=0}^N \left(\frac{h_i}{\sqrt{\tau_i}} I_{n_i}\right) \,
    \mathcal{A}_{\nn}\,
    {\rm diag}_{i=0}^N \left(\frac{h_i}{\sqrt{\tau_i}} I_{n_i}\right)
\]
and, following the same reasoning as in the proof of Theorem \ref{th:general}, we notice that
\[
    \mathcal{A}_{\nn}^{\rm{s}}=
    {\rm diag}_{i=0}^N \left(\frac{h_i}{\sqrt{\tau_i}} I_{n_i}\right)
    D_{\nn} \,\,
    {\rm diag}_{i=0}^N \left(\frac{h_i}{\sqrt{\tau_i}} I_{n_i}\right)
    +R_{\nn}
\]
with rank$(R_{\nn})\le 2\sum_{i=1}^N n_{i,\Gamma}  =  o(n)$, $n=\sum_{i=0}^N n_i$. The statement in \cite[Exercise 5.3]{GS-I} implies again that
\[
    \{\mathcal{A}_{\nn}\}_n^{\rm{s}}
    \quad \mbox{and} \quad
    \left\{
    {\rm diag}_{i=0}^N \left(\frac{h_i}{\sqrt{\tau_i}} I_{n_i}\right)
    D_{\nn} \,\,
    {\rm diag}_{i=0}^N \left(\frac{h_i}{\sqrt{\tau_i}} I_{n_i}\right)
    \right\}_n
\]
share the same distribution, namely
\begin{equation}\label{eq:scaled_matrix_distr}
    \{\mathcal{A}_{\nn}\}_n^{\rm{s}} \sim_\lambda
    \sum_{i=0}^N f^i \psi_{[\hat{r}_{i-1} ,\hat{r}_{i}]},
\end{equation}
with $\hat{r}_{-1}=0$, $\hat{r}_{i}=\hat{r}_{i-1}+r_i$, $\hat{r}_{N}=1$.


In all cases, for symbols $f^i$, where $i=0, \ldots, N$, each $f^i$ maps to $\mathbb{C}^{s_i \times s_i}$, where the dimension $s_i$ is determined by GLT theory as $s_i=(p_i-k_i)^d$. Here, $d$ is the physical dimension of the domain (generally $d=2$ or $d=3$), while $p_i$ is the degree of the polynomial in the Galerkin approach, $k_i$ the global regularity of the approximated solution. More precisely we have $k_i=0$ in standard Lagrangian finite elements, $k_i=-1$ in the discontinuous Galerkin method, $k_i=p_i-1$ when using isogeometric analysis, and $0<k_i<p_i-1$, $p_i\ge 2$, in the case of intermediate regularity: see \cite{garoni2019} for detailed exposition and review paper when $d=1$ and the sections of applications in \cite{block-glt-dD} for problems in multidimensional domains when $d>1$.


\section{Solution strategy and numerical experiments}\label{sec:num}

\subsection{Problem and discretization settings}

For $i=0,\ldots,N$ we set $\sigma_i = 1$ in \eqref{eq::EMI_2}-\eqref{eq::EMI_3} so that $\tau_i=\tau$ for all $i$. If not mentioned otherwise, we set $\tau = 0.01$ in \eqref{eq::EMI_5}-\eqref{eq::f}. We set a passive ionic current in \eqref{eq::f}, i.e. $I_\mathrm{ion}(v) = v$, and an oscillating initial membrane stimulus 
$$v_\mathrm{in}=\frac{1}{2}\sin{(10\|\bm{x}\|^2)} \quad\text{for}\quad \bm{x}\in\Gamma,$$
resulting in the right hand side $g=v_\mathrm{in}(1-\tau)$, with reference to \eqref{eq::f_ij},\eqref{eq::f} and \eqref{eq::f_global}, cf. Figure~\ref{fig:solutions}. We remark that different choices of right hand sides would not alter the essence of this work. In terms of geometry, we consider three different settings: idealized (model A and B) or realistic (model C), described in the following sections. 

\subsection{Model A: idealized nervous system setting}
\begin{figure}
    \centering
    \includegraphics[width=0.32\textwidth]{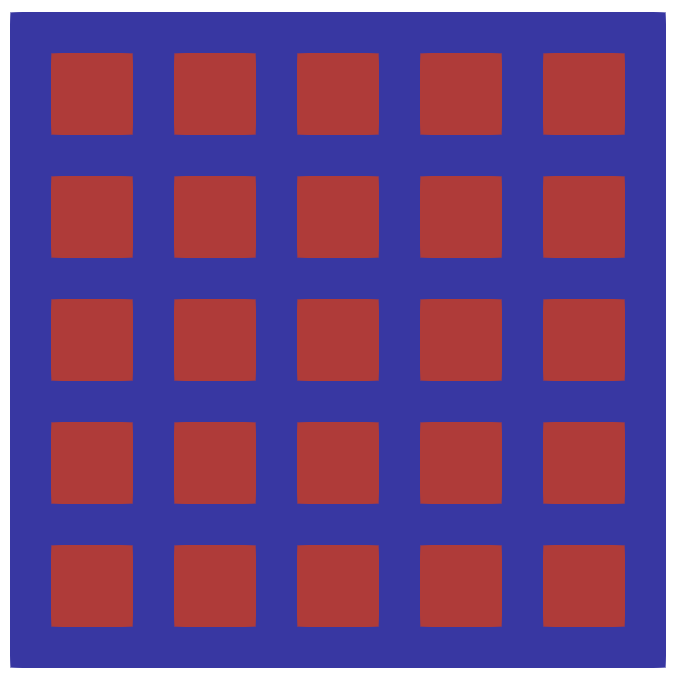}
    \includegraphics[width=0.32\textwidth]{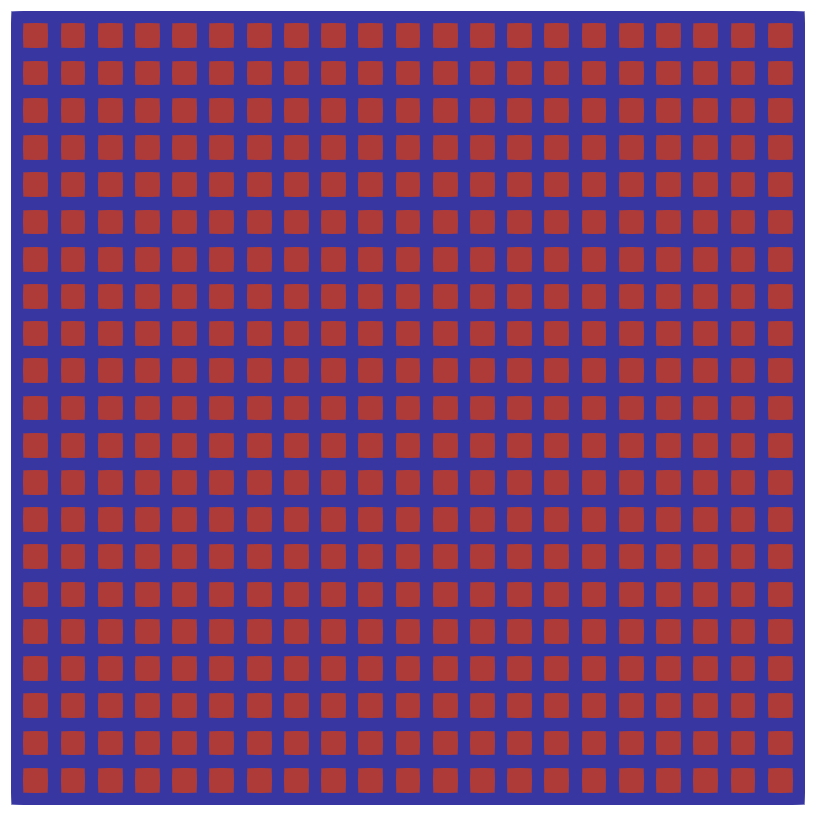}
    \includegraphics[width=0.32\textwidth]{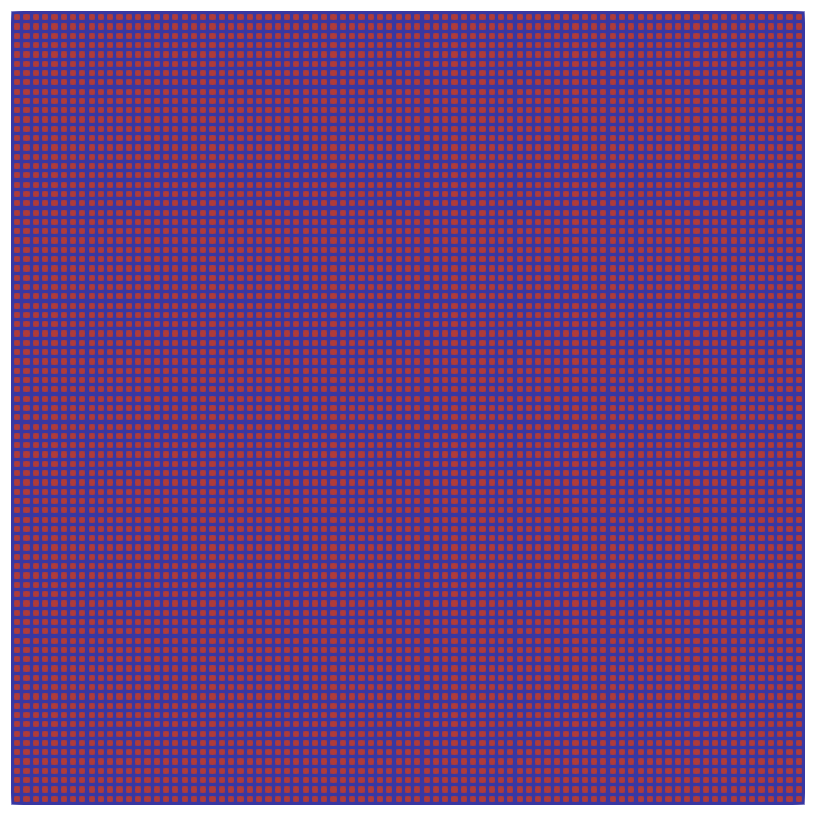}

    $N = 25$ \quad \qquad \qquad\qquad \qquad $N = 441$ \qquad \qquad\qquad  \quad \qquad $N = 7225$
    \caption{Model A: two dimensional geometry varying the number of cells $N$. The extra-cellular space $\Omega_0$ is colored in blue, while cells $\Omega_1,\Omega_2,\ldots,\Omega_N$ in red. No gap junction are present, i.e. $\Gamma_{ij}=varnothing$ for  $i,j\neq 0$.}
    \label{fig:geometry}
\end{figure}

\begin{figure}
    \centering
    \includegraphics[width=0.42\textwidth]{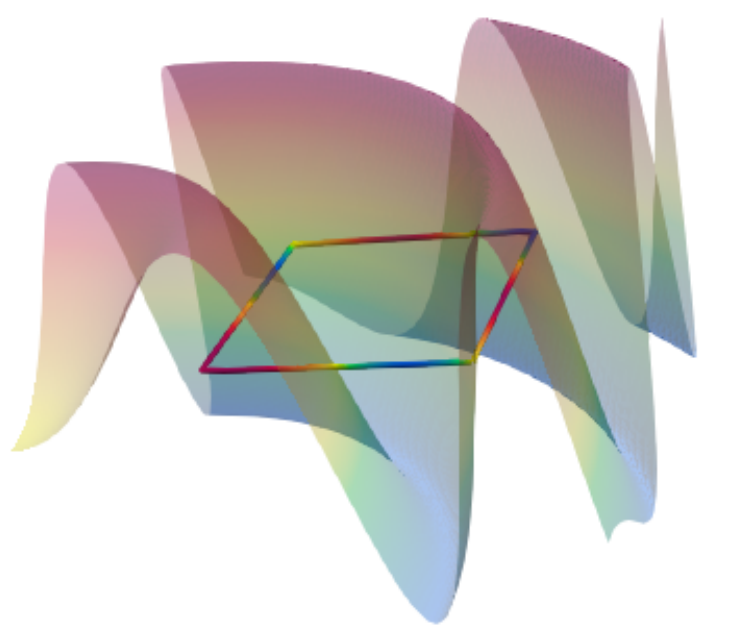}
    \includegraphics[width=0.55\textwidth]{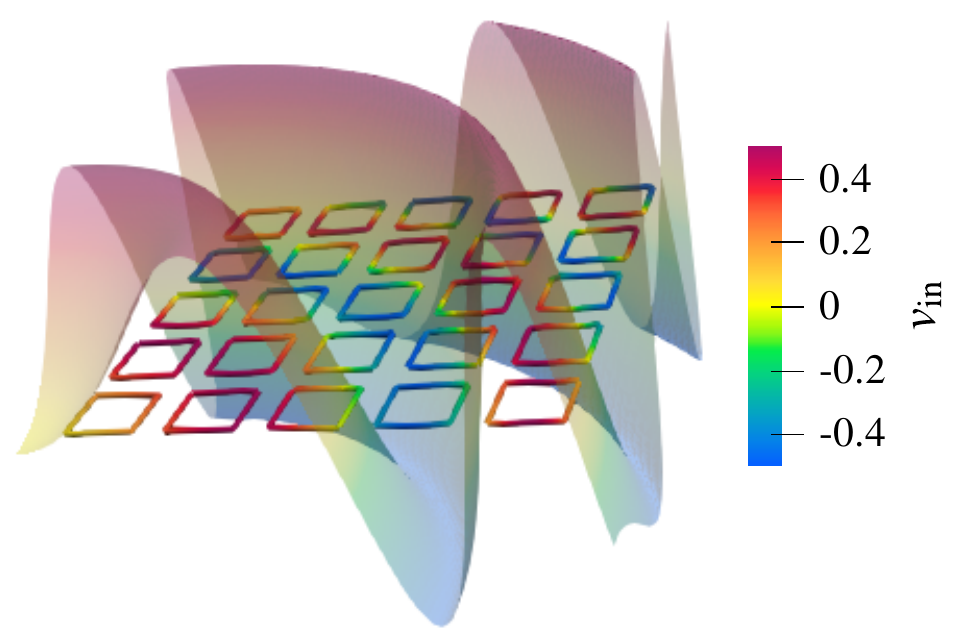}
    \includegraphics[width=0.42\textwidth]{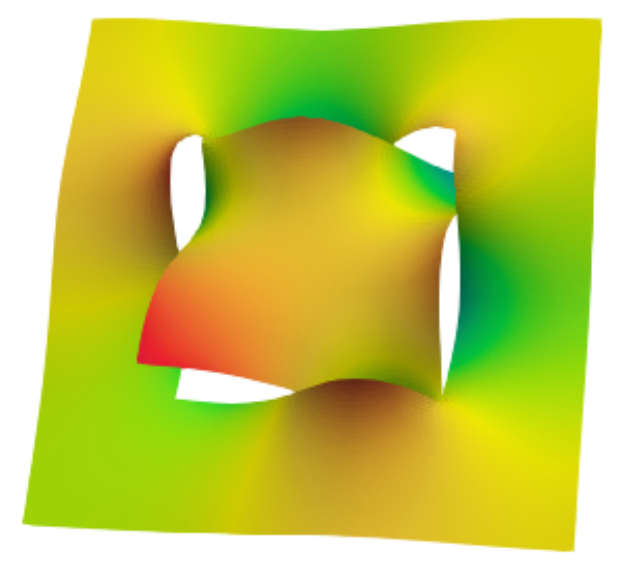}
    \includegraphics[width=0.55\textwidth]{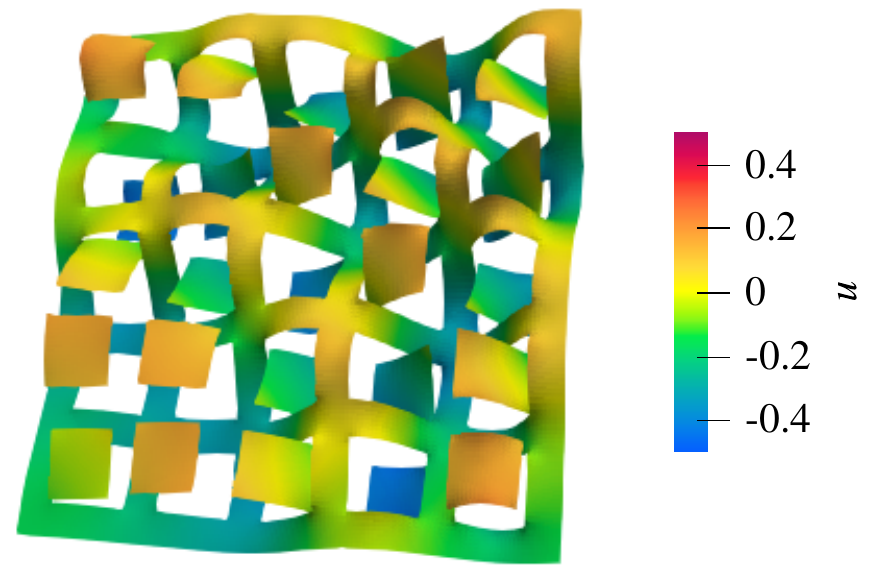}
    $N=1$ \qquad \qquad \qquad \quad \qquad \qquad \qquad $N=25$
    \caption{Model A. First row: overall initial conditions $v_{\mathrm{in}}(x,y)$ for $N=1$ (left) and $N=25$ (right). While $v_{\mathrm{in}}$ is imposed only for $(x,y)\in\Gamma$, we also show it in the whole domain $\Omega$. Second row: corresponding solutions for $N_h = 128$.}
    \label{fig:solutions}
\end{figure}

In the context of Section \ref{sec:nervous_system}, we consider the domain $\Omega = (0,1)^2$.
We  define the $N$-cells partitioning function
$$\psi(x, y) = (x (3N + 1)\mod 3, \,\, y (3N + 1)\mod 3),$$
and the extra- and intracellular regions are respectively given by:
\begin{align*}
  \Omega_0 & = \{(x, y) \in \Omega \mid \min(\psi(x, y)) \leq 1\},  \\
  \Omega_1\cup\Omega_2\cup\ldots \cup\Omega_N & = \{(x, y) \in \Omega \mid \min(\psi(x, y)) \geq 1\}
\end{align*}
cf. Figure~\ref{fig:geometry}.
The domain $\Omega$ is discretized with with a triangular tessellation with $N_h$ elements per side, with a total number of degrees of freedom $n=n_0+n_\mathrm{in}=\sum_{i=0}^Nn_i$ given $(N_h+1)^2$ grid points. Assuming that $N_h$ is a power of 2, the number of cells is chosen of the form $N=\left(\frac{2^{2k}-1}{3}\right)^2$, with $k = 1,2,\dots , \lceil(\log_2 N_h)/2\rceil$, so that the geometry is compatible with the tessellation and the interfaces $\Gamma_i$ are well defined on the mesh. In this setting, the same number of degrees of freedom $n_i = n_\mathrm{in}/N$ corresponds to all cells $i=1,\ldots,N$, and $n_0, n_\mathrm{in}$, and $n_\Gamma$ depend on $N$ cf. Table~\ref{tab:geometry} and Figure~\ref{fig:geometry}. The overall solution, denoted by $u$, including $u_0,\ldots,u_N$ and the membrane initial condition $v_{\mathrm{in}}$ are shown in Figure~\ref{fig:solutions} for $N=1$ and $N=25$.

\begin{table}[]
    \centering
    \begin{tabular}{l|r|r|r|r|r}
        Cells $N$ & 1 & 25 & 441 & 7225 & 116281 \\
         \hline
         Extracell. dofs $n_0$  & 789504  & 647400 & 626824 & 696600 & 934344  \\
         Intracell. dofs $n_\mathrm{in}$ & 263169 & 416025 & 480249 & 924800 & 1046529  \\
         Membrane dofs $n_\Gamma$ & 2048 & 12800 & 56448 & 231200 & 930248 \\
         Total dofs $n=n_0+n_\mathrm{in}$ & 1052673 & 1063425 &  1107073 & 1621400 & 1980873 \\
         $n_\Gamma/n$ & 0.002 &  0.012  & 0.051 & 0.143 & 0.470
    \end{tabular}
    \caption{Model A: example of number of degrees of freedom corresponding to various regions for $N_h=1024$, i.e. $(N_h+1)^2 = 1050625$ grid points, varying the number of cells $N$. }
    \label{tab:geometry}
\end{table}

\subsection{Model B: idealized myocytes}
\begin{figure}
    {\centering
    \includegraphics[width=0.33\textwidth]{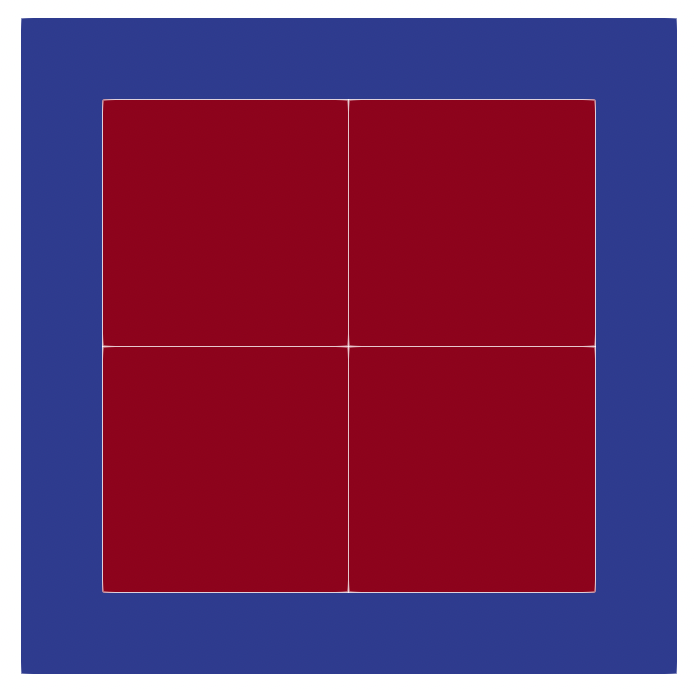}
    \qquad
    \includegraphics[width=0.33\textwidth]{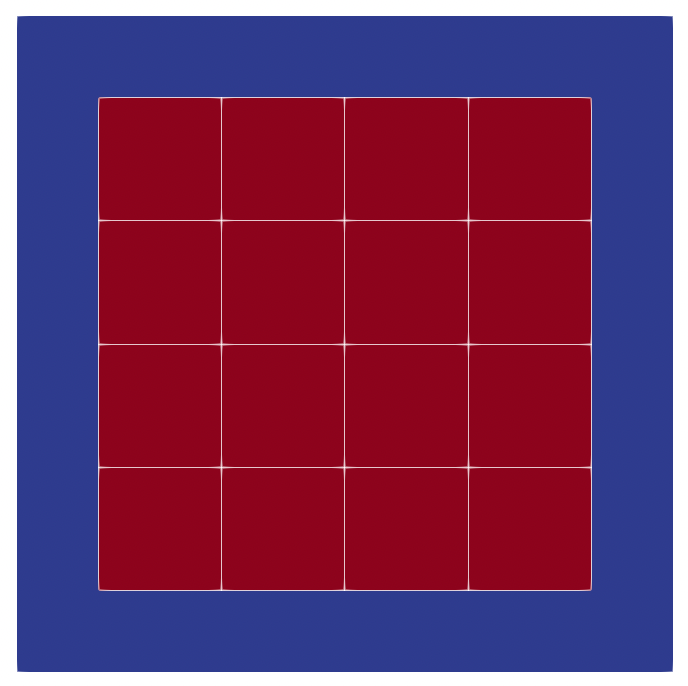}

    }
    
    \qquad  \includegraphics[width=0.4\textwidth]{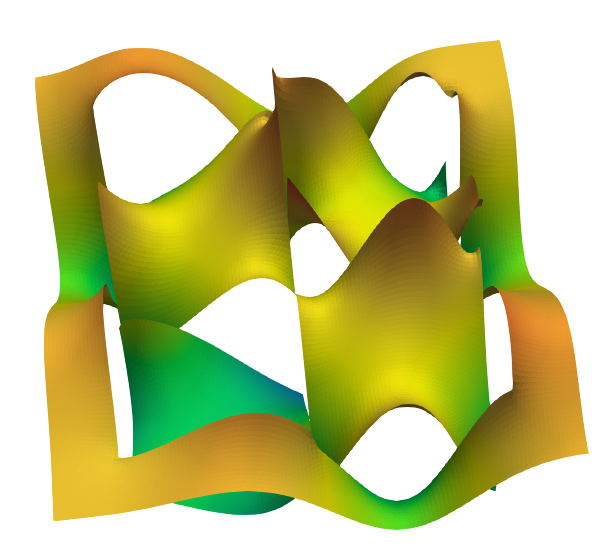}
    \qquad
    \includegraphics[width=0.57\textwidth]{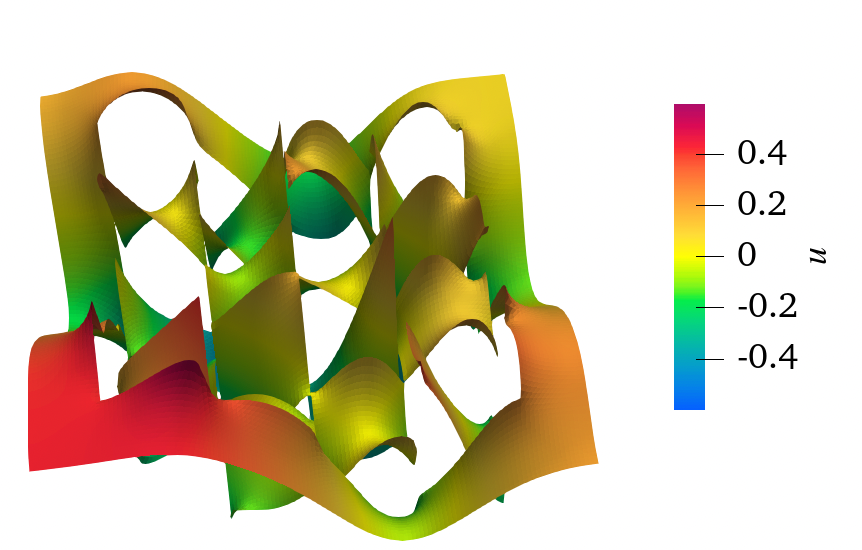}
    
    \quad \quad \quad \qquad \qquad \qquad $N = 4$ 
    \quad \quad \quad \quad \quad \quad \quad \quad \qquad \qquad \qquad $N = 16$ 
    \caption{Model B. First row: 2D geometry with $N=4$ and $N=16$ cells; the extra-cellular space $\Omega_0$ is colored in blue, while cells $\Omega_1,\ldots,\Omega_N$ in red and membranes in white. Second row: corresponding solutions $u_0,\ldots,u_N$ for $N_h = 512$.}
    \label{fig:geometry2}
\end{figure}

In the context of Section~\ref{sec:cardiac_setting}, we consider the domain $\Omega = (0,1)^2$ and the intracellular subdomains $\Omega_1\cup\Omega_2\cup\cdots\cup\Omega_N = (\sfrac{1}{8},\sfrac{7}{8})^2$.
Again, $\Omega$ is discretized using a triangular tessellation with $N_h$ elements per side, with $N_h$ a power of two. The $N$ cells are disposed in a regular $\sqrt{N}\times\sqrt{N}$ grid. Each cell corresponds to $(1+\frac{3N_h}{4\sqrt{N}})^2$ degrees of freedom, so that the total number of intracellular degrees of freedom is $n_\text{in}=N\cdot(1+\frac{3N_h}{4\sqrt{N}})^2$ and extracellular ones $n_0=\frac{7}{2}N_h(N_h/8+1)$. Examples of this geometry and corresponding solutions are reported in Figure~\ref{fig:geometry2}.

\subsection{Model C: mouse visual cortex tissue}
\begin{figure}
    \centering
    \includegraphics[width=0.4\textwidth]{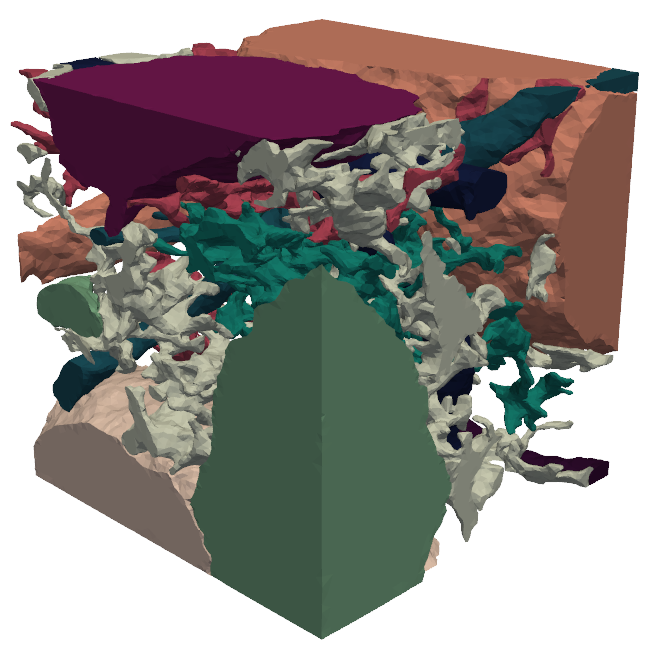} \qquad
    \includegraphics[width=0.4\textwidth]{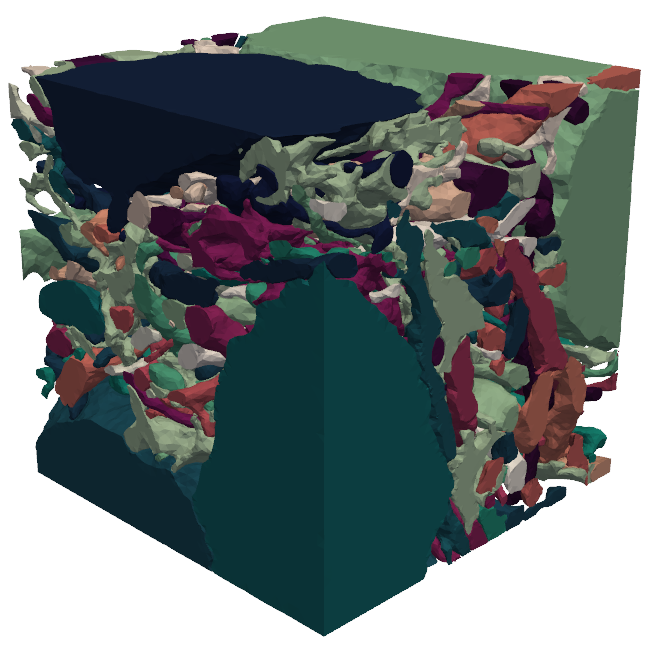}\\
     $N=10$ \qquad \qquad \qquad \qquad \qquad \qquad \qquad \qquad \qquad $N=50$\\
     \vspace{0.7cm}
    \includegraphics[width=0.4\textwidth]{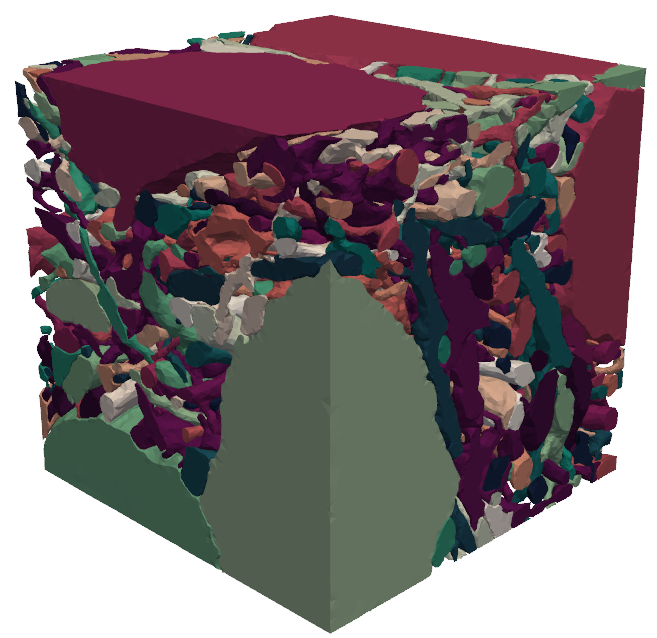} \qquad
    \includegraphics[width=0.4\textwidth]{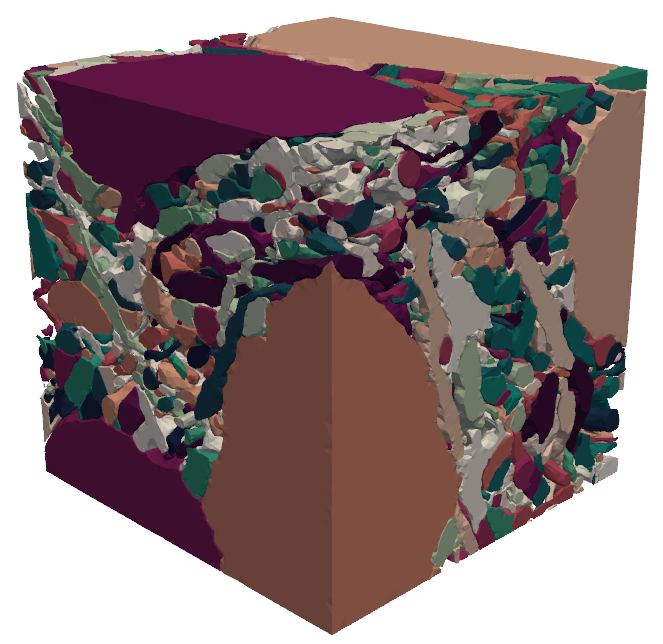}\\
    $N=100$ \qquad \qquad \qquad \qquad \qquad \qquad \quad \qquad \qquad $N=200$
    \caption{Model C. For $N=10$ the largest ten cells in the dataset are included; four neuronal somas and few glial cells are visible. As $N$ increases the extracellular space is filled with more and more glial cells.}
    \label{fig:model_C_geo}
\end{figure}

\begin{figure}
    \centering
    \includegraphics[width=0.49\linewidth]{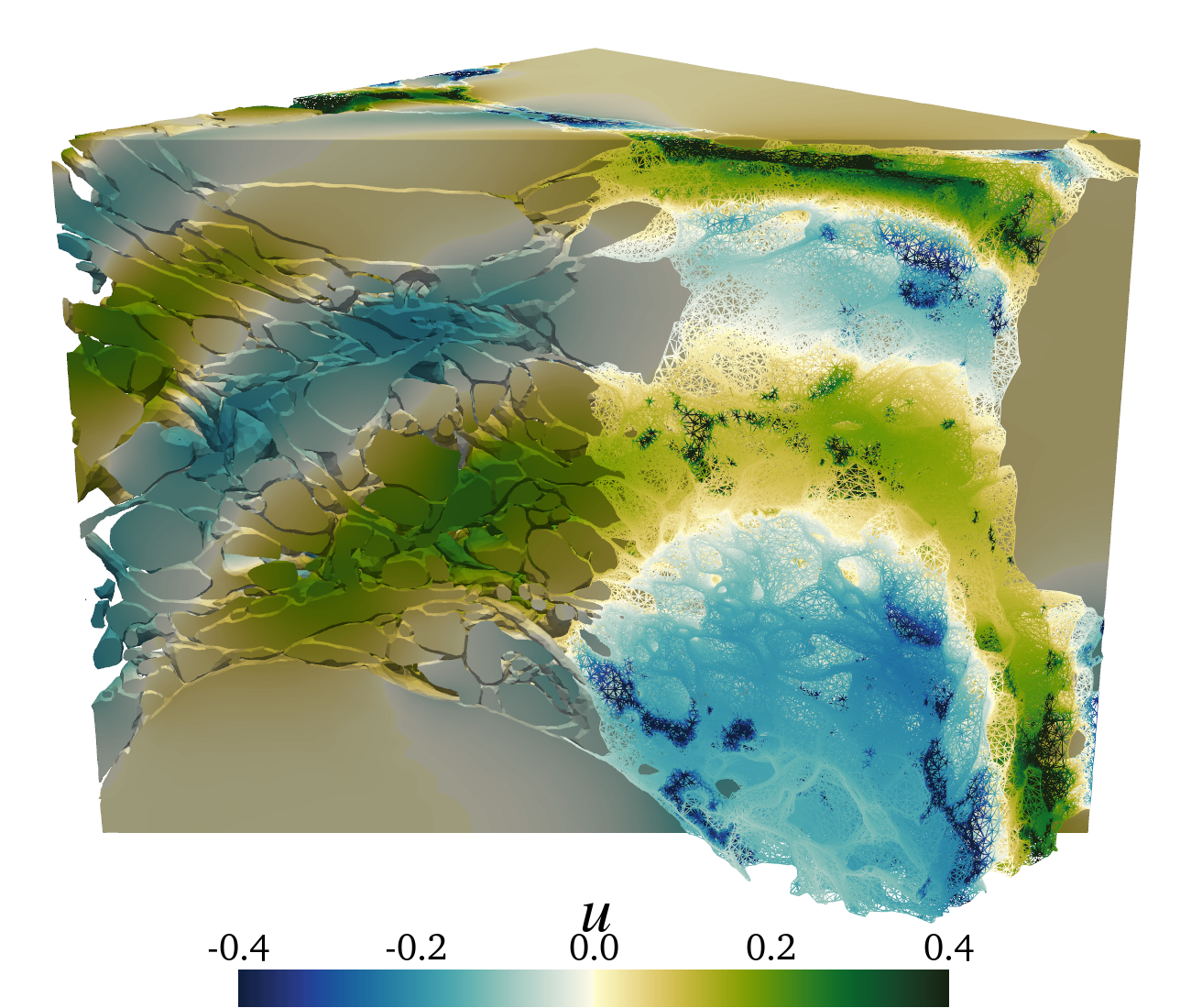}
    \includegraphics[width=0.49\textwidth]{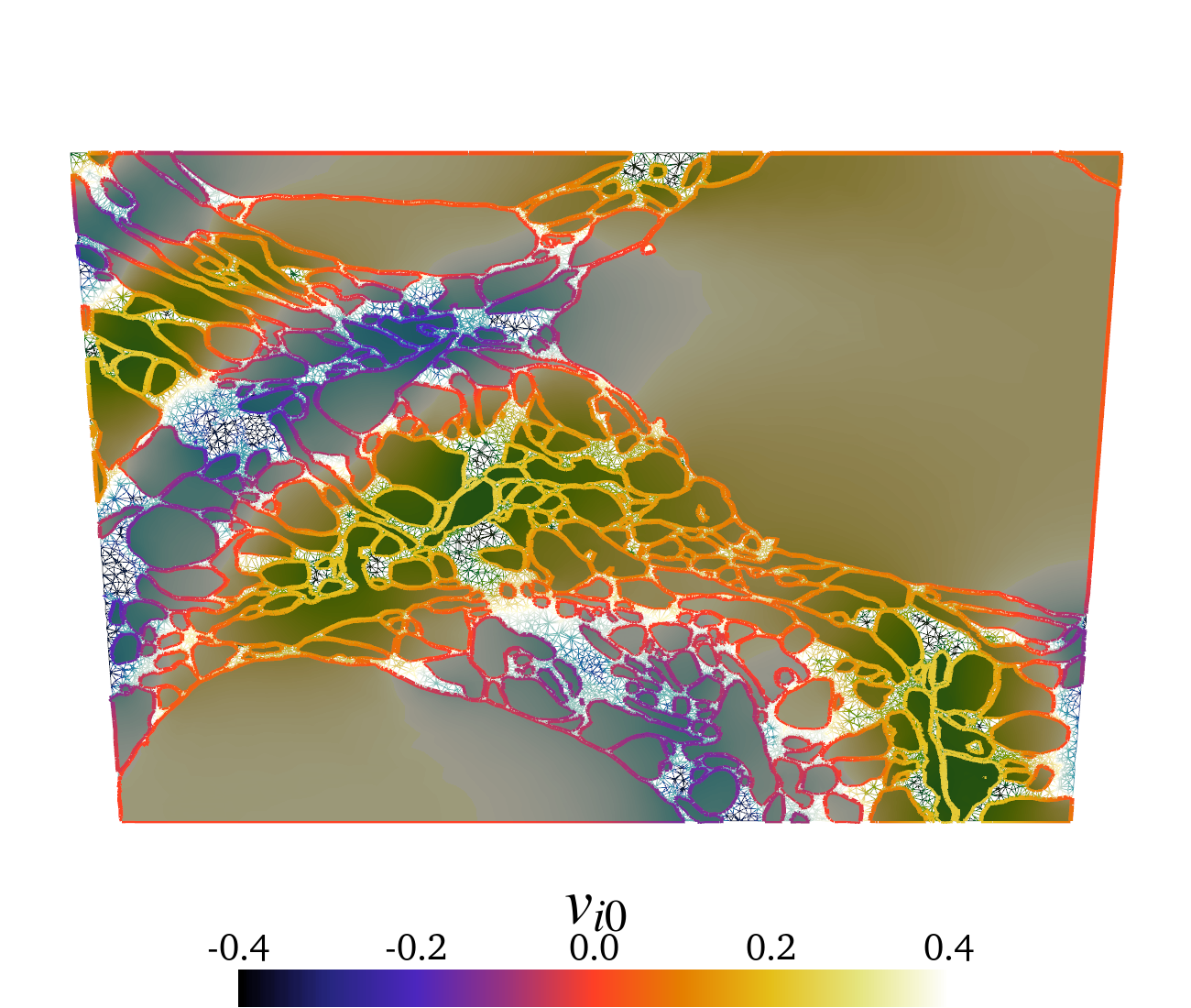}\\
    \caption{Model C with $N=200$ cells: potentials $u_0,\ldots,u_N$ on the clipped 3D domain (left) and a 2D domain slice (right), additionally highlighting the membrane potential $v_{i0}$ with $i=1,\ldots,N$. For the extracellular space $\Omega_0$, the corresponding mesh is shown. The colorbar for the potential applies to the intra- and extracellular potential in both plots.}
    \label{fig:model_C}
\end{figure}

We consider tetrahedral meshes of a dense reconstruction of the mouse visual cortex at extreme resolution. The dataset is based on the Cortical MM$^3$ dataset\footnote{\url{https://www.microns-explorer.org/cortical-mm3}} centred at position (225182,107314,22000) with a resolution of $32\times32\times40\,\text{nm}^3$. The meshed domains are cubes with side length 10 $\mu$m and include the largest 10, 50, 100 and 200 cells in the respective tissue volume, respectively, cf. Figure~\ref{fig:model_C_geo} and Table~\ref{tab:geometry_C} for geometrical data. The domain is scaled so that $\Omega=(0,1)^3$. We show a model C solution and membrane potential in Figure~\ref{fig:model_C}. We remark that the geometrical setting is different w.r.t. Table~\ref{tab:geometry} for model A, since the tessellation in model C is changing to accommodate for new cells and the insertion of new cells does not modify the existing ones.

\begin{table}[]
    \centering
    \begin{tabular}{l|r|r|r|r|r}
        Cells $N$ & 10 & 50 & 100 & 200\\
         \hline
         Extracell. dofs $n_0$ & 185918 & 371684 & 467421 & 635518 \\
         Intracell. dofs $n_\mathrm{in}$ & 139677 & 364293 & 487399 & 319302\\
         Membrane dofs $n_\Gamma$ & 104547 & 298463 & 423930 & 609369\\
         Total dofs $n=n_0+n_\mathrm{in}$ & 325595 &  735977 & 954820 & 1327301 \\
         $n_\Gamma/n$ & 0.32 & 0.41 & 0.44 & 0.46
    \end{tabular}
    \caption{Model C: number of degrees of freedom corresponding to various regions varying the number of cells $N$.}
    \label{tab:geometry_C}
\end{table}

\subsection{Implementation and solution strategies}
We use FEniCS \cite{alnaes2015fenics,logg2012automated} for parallel finite element assembly and \texttt{multiphenics}\footnote{\url{https://multiphenics.github.io/index.html}} to handle multiple meshes with common interfaces and the corresponding mixed-dimensional weak forms\footnote{For the corresponding open software visit \url{https://github.com/pietrobe/EMIx}}. FEniCS wraps PETSc parallel solution strategies. For multilevel preconditioning, we use a single iteration of \texttt{hypre boomerAMG} algebraic multigrid\footnote{We use default options, only for $d=3$ we set the strong coupling threshold to 0.5.} (AMG$_1$) \cite{falgout2002hypre}. 
For compara\-tive studies, we use incomplete LU (ILU) preconditioning with zero fill-in. Parallel ILU is obtained via a block Jacobi preconditioner, where each block inverse is approximated by ILU. For iterative strategies, as stopping criteria, we use a tolerance of $10^{-9}$ for the relative residual. Given the pure Neumann boundary conditions \eqref{eq:bc}, uniqueness is enforced by setting a point-wise boundary condition $u_0(\mathbf{0})=0$.
Run-times are reported for a 2021 MacBook Pro with an M1 Pro chip, 8 cores, and 16 GB of memory. When not mentioned otherwise, experiments are performed in serial.

Following \eqref{eq::SMW}, the analysis in \cite{benedusi2023EMI}, and Theorem~\ref{th:general}, we also introduce the block diagonal preconditioner, for $\epsilon>0$,
\begin{equation}\label{eq::P}
P_{\epsilon,\nn} = \tau	
 \begin{bmatrix}
 A_0 + \epsilon\widetilde{M}_0 &  &  &  & \\
	&  A_1 + \epsilon\widetilde{M}_1 &  &  &   \\
     &   & A_2 + \epsilon\widetilde{M}_2 &  &   \\
     &   &  & \ddots &  \\
     &   & &  & A_N + \epsilon\widetilde{M}_N  \\	
	\end{bmatrix}\in\mathbb{R}^{n\times n},
\end{equation}
where the bulk Laplacians $A_i$, for $i=0,\ldots,N$, are regularized adding the corresponding bulk mass matrices
\begin{equation*}
    \widetilde{M}_i=\left[\int_{\Omega_i}\phi_{i,\ell}(\xx)\phi_{i,k}(\xx)\,\mathrm{d}\xx\right]_{\ell,k=1}^{n_i}\in\mathbb{R}^{n_i\times n_i}.
\end{equation*}
Let us remark that, while $M_i$ corresponds to degrees of freedom in $\Gamma_i$, the matrix $\widetilde{M}_i$ refers to $\Omega_i$, and is necessary as a full-rank correction to $A_i$, which is rank-1 deficient, lacking boundary conditions. 

From a theoretical point of view, a few facts have to be emphasized as a consequen\-ce of the theoretical analysis.
First, we notice that 
\begin{equation}\label{eq::A-P}
\{\mathcal{A}_{\nn}-P_{\epsilon,\nn}\}_n \sim_\lambda 0 
\end{equation}
for any $\epsilon>0$ thanks to the different scaling of the mass and stiffness matrices, in accordance with the analysis in Section \ref{ssec:symbol} and especially in Section \ref{ssec:symbol-bis}. Therefore the two matrix-sequences $\{\mathcal{A}_{\nn}\}_n,\ \{P_{\epsilon,\nn}\}_n$ share the same distribution and more importantly
\begin{equation}\label{eq::invP*A}
\{P_{\epsilon,\nn}^{-1}\mathcal{A}_{\nn}-I_n\}_n \sim_\lambda 0 
\end{equation}
because both $\{\mathcal{A}_{\nn}\}_n,\ \{P_{\epsilon,\nn}\}_n$ are sparsely vanishing (see \cite[Definition 8.2]{GS-I}), i.e. its common distribution has a determinant vanishing in a set of zero Lebesgue measure in agreement with \cite[Proposition 2.27]{block-glt-dD}.
We remark that (\ref{eq::invP*A}) is an important indication that the considered preconditioning strategy is effective. For more analysis, also in connection with outliers, we refer to \cite[Section 4.4]{benedusi2023EMI}, where the case of a single cell ($N=1$) was treated in detail.

\subsection{Results}
We test the convergence of various preconditioners for the CG solver, to highlight which is the more appropriate, depending on the application of interest. The block diagonal preconditioner $P_{\epsilon,\nn}$ is solved to full precision with Cholesky from MUMPS and labelled as $P_{\epsilon}$-CG, using $\epsilon=10^{-4}$ for models A,B, and $\epsilon=1$ for model C. As a preliminary experiment, in Table~\ref{tab:results1}, we consider both model A and B, with a fixed number of cells $N$ and refining the mesh. Compared to CG or ILU-CG, both the block preconditioned CG and the multilevel preconditioner are essentially robust w.r.t. $N_h$.  Interestingly, $P_{\epsilon}$-CG is less effective and robust in the cardiac setting (Model B), for which the {latter theory cannot be directly applied}. The direct computation of $P_\epsilon^{-1}\mathcal{A}_n$ spectrum shows a constant number of outliers (e.g. outside the interval $[1-\epsilon,1+\epsilon]$), which are not captured by the spectral analysis, and can explain the observed convergence. The analysis of the outliers is quite involved especially in the preconditioned case and need specific tools; see \cite{prec-out1,prec-out2,SeTi} and references therein. In Table~\ref{tab:results1.1}, we asses robustness w.r.t. the temporal discretization parameter $\tau$ for the preconditioning strategies. As expected from equation~\eqref{Di-expression} and the discussion in Section~\ref{ssec:symbol-bis}, the block diagonal preconditioner is not robust as $\tau\to 0$. Surprisingly, AMG$_1$-CG is effective also in this scenario.

\begin{table}[ht]
    \centering
    \begin{tabular}{l|l|r|r|r|r|r}
         & $N_h$ & 64 & 128 & 256 & 512 & 1024 \\
         \hline
         \hline
        Model A & CG & 392 & 743 & 1437 & 2725 & 5471 \\
         $N=441$ & ILU-CG & 103 & 178 & 320 & 600 & 1694 \\         
        & $P_\epsilon$-CG & 46 & 53 & 58 & 61 & 61 \\
        & AMG$_1$-CG  & 8 & 8 & 9 & 10 & 10 \\
         \hline
          \hline
        Model B & CG & 535 & 954 & 1810 & 3506 & 6931 \\  
          $N=576$& ILU-CG & 129 & 188 & 327 & 627 & 1235 \\ 
         & $P_\epsilon$-CG & 1040 & 1129 & 1221 & 1256 & 1352  \\
        & AMG$_1$-CG  & 9 & 9 & 10 & 10 & 12         

    \end{tabular}
    \caption{Iterations to convergence for various iterative methods with a fixed number of cells $N=441$ (resp. $N=576$) for model A (resp. B) and refining spatial discretization via $N_h$.}
    \label{tab:results1}
\end{table}

\begin{table}[ht]
    \centering
    \begin{tabular}{l|r|r|r|r|r}
         $\tau$ & 0.1 & 0.01 & 0.001 & 0.0001 & 0.00001 \\
         \hline
         ILU-CG  & 1208 & 943 & 786 & 763 & 787 \\
         $P_\epsilon$-CG & 27 & 61 & 139 & 211 & 226 \\
         AMG$_1$-CG & 11 & 9 & 8 & 8 & 7

    \end{tabular}
    \caption{Model A: iterations to convergence for various iterative methods with fixed spatial discretization $N_h=512$, the number of cells $N=441$, and varying time step $\tau$. Results using 8 cores.}
    \label{tab:results1.1}
\end{table}
In the following results, we consider an increasing number of cells $N$ for a fixed discretization. These experiments are of interest for at least two reasons:
\begin{itemize}
    \item[(i)] relevant portions of realistic tissue are typically densely populated by many different cells. In particular, in the neuroscience context, membranes have complex morphologies and the membrane degrees of freedom $n_\Gamma$ are not negligible w.r.t. $n$, unless aggressive mesh refinement is used, cf. Figure~\ref{fig:model_C_geo} and Table~\ref{tab:geometry_C};
    \item[(ii)] as $N$ increases, also $n_\Gamma$ increases, challenging the hypothesis of Theorem~\ref{th:general}.
\end{itemize}
According to Tables~\ref{tab:results2}-\ref{tab:results3}, the block diagonal preconditioner $P_{\epsilon}$ is effective only for model A, and the monolithic multilevel strategy AMG$_1$-CG is effective and robust in both cases.This superior performance of AMG$_1$-CG can be attributed to the spectral properties of the whole coefficient matrix-sequence, which exhibits Laplacian-like characteristics in the related eigenvalue distribution, regardless of the specific setting. As expected from theory, $P_{\epsilon}$-CG is less effective as $N$ increases, since the first hypothesis in Theorem~\ref{th:general} is no longer satisfied. 
Furthermore, as already observed at the beginning of Section \ref{sec:cardiac_setting} and in Remark \ref{rem:cardiac-trid}, the block diagonal preconditioner cannot be satisfactory for model B, since it does not capture the spectral distribution of the original coefficient matrices.

Similarly, for the realistic model C, we report run-times and iteration counts in Table~\ref{tab:results3} for an increasing number of cells $N$, showing robustness and efficiency for both the ILU and the AMG$_1$ preconditioners (competitive, e.g., w.r.t. the timing results in \cite{jaeger2021efficient}). In this case the block preconditioner $P_{\epsilon}$ is never viable. 
Let us remark that multilevel monolithic solution strategies showcase good scalability when employed on large-scale problems using massively parallel machines \cite{benedusi2024scalable,benedusi2021fast}.

\begin{table}[ht]
    \centering
    \begin{tabular}{l|r|r|r|r|r}
         $N$ & 1 & 25 & 441 & 7225 & 116281 \\
         \hline         
         ILU-CG & 1796 & 1713 & 1920 & 1909 & 2353\\
        $P_\epsilon$-CG & 57 & 129 & 276 & 610 & 1112  \\
         AMG$_1$-CG & 9 & 9 & 11 & 11 & 8        
    \end{tabular}
    \caption{Model A: iterations to convergence of various iterative schemes fixing the spatial discretization with $N_h=1024$, i.e. $(N_h+1)^2 = 1050625$ grid points. Results using 8 cores.}
    \label{tab:results2}
\end{table}

\begin{table}[ht]
    \centering
    \begin{tabular}{l|r|r|r|r|r}
         $N$ & 1 & 16 & 256 & 576 & 4096\\
         \hline         
         ILU-CG          & 577 & 588 & 614 & 627 & 704\\
         $P_\epsilon$-CG & 65 & 305 & 969 & 1256 & 2295 \\
         AMG$_1$-CG      & 8 & 9 & 10 & 10 & 11 \\
    \end{tabular}
    \caption{Model B: iterations to convergence of various iterative schemes fixing the spatial discretization with $N_h=512$, i.e. $(N_h+1)^2 = 263169$ grid points.}
    \label{tab:results3}
\end{table}

\begin{table}[ht]
    \centering
    \begin{tabular}{l|r|r|r|r}
         $N$ & 10 & 50 & 100 & 200\\
         \hline   
         ILU-CG          & 0.9 [248] & 3.0 [320] & 3.9 [257] & 6.3 [368] \\
         $P_\epsilon$-CG & 66 [1069] & 397 [2220]  & 812 [2344] & 1224 [2452]\\
         AMG$_1$-CG      & 1.4 [11] & 4.2 [12] & 6.3 [13] & 9.2 [12]
    \end{tabular}
    \caption{Model C: run-times (s) and iterations (in square brackets) to convergence of ILU and AMG$_1$ preconditioners, compered to a direct approach, increasing the number of cells $N$ with geometries according to Table~\ref{tab:geometry_C}. Results using 8 cores.}
    \label{tab:results4}
\end{table}

\section{Concluding remarks}\label{sec:final}
We described a numerical scheme for approximating the solution of the EMI equations under the assumptions of $N$ distinct cells. We showed the
structure and spectral features of the coefficient matrices obtained from Galerkin discretizations, which are spectrally described by an appropriate linear combi\-nation of Laplacian symbols, cf. Theorem~\ref{th:general}. This information suggests the use of both block diagonal preconditioners and multilevel strategies to
accelerate Krylov solvers, when the spectral analysis hypotheses are satisfied. Realistic settings might challenge these hypotheses
and we can expect degradation of block preconditioner performance, as we observe in the numerical experiments. In practice, a single monolithic multigrid iteration is the most robust strategy, according to the presented tests. Geometrical dependencies on convergence have yet to be fully understood. In general, many open problems remain, e.g., the precise study of the outliers, of their number, and of the asymptotic behavior with respect to the matrix-size and fineness parameters: using, e.g., the tools given in \cite{prec-out1,prec-out2,SeTi}. Finally, with regard to the SMW formulae (\ref{eq::SMW-eps}) and (\ref{eq::SMW}), beside their algebraic structure, it is of interest to perform a GLT based spectral analysis. The difficulty is represented by the presence of rectangular terms (as  ${U}_{\nn}, {V}_{\nn}, \widetilde{U}_{\nn}, \widetilde{V}_{\nn}$), but this can be overcome by using quite recent extensions of the GLT theory as that in \cite{rect-glt}.
However, the GLT analysis and the design of symbol-based preconditioners as in \cite{NS1-elonged} are not immediate and they are beyond the scope of the current work.


\subsection*{Acknowledgments}
 Pietro Benedusi and Marius Causemann acknowledge sup\-port from the Research Council of Norway via FRIPRO grant \#324239 (EMIx). This work is also financially supported by the Center for Computational Medicine in Cardiology. Paola Ferrari and Stefano Serra-Capizzano are partially supported by the Italian Agency INdAM-GNCS. Furthermore, the work of Stefano Serra-Capizzano is funded from the European High-Performance Computing Joint Undertaking  (JU) under grant agreement No 955701. The JU receives support from the European Union’s Horizon 2020 research and innovation programme and Belgium, France, Ger\-many, and Switzerland.
  Stefano Serra-Capizzano is also grateful for the support of the Laboratory of Theory, Economics and Systems – Department of Computer Science at Athens University of Economics and Business. The work of Paola Ferrari was carried out within the framework of the project HubLife Science – Digital Health (LSH-DH) PNC-E3-2022-23683267 - Progetto DHEAL-COM – CUP:D33C22001980001, founded by Ministero della Salute within ``Piano Nazionale Complementare al PNRR Ecosistema Innovativo della Salute - Codice univoco investimento: PNC-E.3''. Special thanks are also extended to J\o rgen Dokken, Halvor Herlyng, Marie E. Rognes, Patrick Zulian, and Hardik Kothari for the valuable discussions.

\section*{Data availability}
The used open software is available at \url{https://doi.org/10.5281/zenodo.10728462} \cite{benedusi2024knpemi-zenodo}.
Data sets generated during the current study are available from the corresponding author on reasonable request.

\section*{Conflict of interest}
The authors declare that they have no conflict of interest.

 \bibliographystyle{siamplain}
\bibliography{references}

\end{document}